\documentclass[a4paper]{article}

\usepackage{a4wide}
\usepackage{latexsym}
\usepackage[UKenglish]{babel}
\usepackage[colorlinks=true,linkcolor=blue,citecolor=blue]{hyperref}
\usepackage{amsfonts,amsmath,amsthm,mathtools,bbm,cool}
\usepackage{soul} 

\usepackage{lipsum}
\usepackage{amsfonts}
\usepackage{graphicx,subfigure,caption}
\usepackage{epstopdf}
\graphicspath{{./}{figures/}}
\usepackage{algorithmic}
\ifpdf
\DeclareGraphicsExtensions{.eps,.pdf,.png,.jpg}
\else
\DeclareGraphicsExtensions{.eps}
\fi
\usepackage[UKenglish]{babel}
\usepackage{tikz}
\usepackage{pgfplots}
\usepackage{hhline}
\usepackage{tabularx} 
\usepackage{tablefootnote}
\usepackage{multirow}

\usepackage{makecell}
\usepackage{colortbl}
\usepackage{accents}
\newcommand{\dbtilde}[1]{\accentset{\approx}{#1}} 

\usepackage{tikz}
\usepackage{pgfplots}

\usepackage{hhline}
\usepackage{tabularx} 
\usepackage{multirow}

\usepackage{makecell}

\newtheorem{theorem}{Theorem}[section]
\newtheorem{proposition}[theorem]{Proposition}

\newtheorem{lemma}{Lemma}[section] 

\theoremstyle{remark}\newtheorem{remark}[theorem]{Remark}
\newtheorem{assumption}{Assumption}[section]

\usepackage{todonotes}

\newcommand{\tn}{\textnormal}
\newcommand{\rev}{\textcolor{black}}

\newcommand{\numberset}{\mathbb}
\newcommand{\EE}{\numberset{E}} 
\newcommand{\PP}{\numberset{P}}

\newcommand{\B}{\mathcal{B}}
\newcommand{\nuy}{\nu^\mathbf{Y}}

\newcommand{\A}{\mathcal{A}}
\newcommand{\C}{\mathcal{C}}

\newcommand{\N}{\mathcal{N}}

\newcommand{\F}{\mathcal{F}}

\newcommand{\pd}{\partial}

\newcommand{\R}{\mathbb{R}}
\newcommand{\Rd}{\mathbb{R}^d}
\newcommand{\Rk}{\mathbb{R}^k}

\usepackage{mathrsfs}

\newcommand\restr[2]{{
		\left.\kern-\nulldelimiterspace 
		#1 
		\littletaller 
		\right|_{#2} 
}}
\newcommand{\littletaller}{\mathchoice{\vphantom{\big|}}{}{}{}}

\newcommand{\xs}{x_{\star}}

\newcommand{\bfy}{\mathbf{y}}
\newcommand{\by}{\overrightarrow{\mathbf{y}}}

\newcommand{\bfY}{\mathbf{Y}}
\newcommand{\bY}{\overrightarrow{\mathbf{Y}}}

\begin{document}
\title{Consensus-based algorithms for stochastic optimization problems}

\author{
	Sabrina Bonandin \\
	{\small	Institute for Geometry and Practical Mathematics},
		{\small RWTH Aachen University, Germany} \\
		{\small\tt bonandin@eddy.rwth-aachen.de} \\
	Michael Herty \\
		{\small	Institute for Geometry and Practical Mathematics},
		{\small RWTH Aachen University, Germany} \\
		{\small	Extraordinary Professor, Department of Mathematics and Applied Mathematics},\\
		{\small University of Pretoria, South Africa}\\
		{\small\tt herty@igpm.rwth-aachen.de} 
		}  
\date{20th November 2025}

\maketitle

\begin{abstract}
We address an optimization problem where the cost function is the expectation of a random mapping.
To tackle the problem two approaches based on the approximation of the objective function by consensus-based particle optimization methods on the search space are developed. 
The resulting methods are mathematically analyzed using a mean-field approximation, and their connection is established. 
Several numerical experiments show the validity of the proposed algorithms and investigate their rates of convergence.
\medskip


 \end{abstract}


\section{Introduction}
\label{sec:intro}
The interest in addressing optimization problems that accommodate for uncertainties has risen in the last decade \cite{chen2011stochastic}. 
When considering optimization problems under uncertain information, there are mainly two choices to be made \cite{bianchi2009survey}: first, the way to formalize the variability and, second, the time the uncertain information is revealed with respect to the decision-making. 
In this manuscript we restrict to the class of static stochastic optimization problems (sSOPs \footnote{The definition of sSOPs is not unanimous. We refer to \cite{bianchi2009survey}.}), namely to minimization tasks where the uncertain information is described by means of random variables and where the optimization effort is applied before the random event has taken place. More precisely, we consider problems of the type:
\begin{equation}
	\label{eqi: min problem main}
	\min_{x \in \Rd} \{ f(x) := \EE_\PP[F(x,\mathbf{Y})]\}.
\end{equation}
Here $d\ge 1$, $\mathbf{Y}$ is a random variable defined on the probability space $ (\Omega, \mathcal{A}, \PP) $ and taking values in a set $ E \subset \Rk, k \ge 1$. Denoting by $\B(E)$ the Borel set of $E$, we call $\nu^{\mathbf{Y}}: \B(E) \to [0,1]$ its law.
Further, $F:\Rd \times E \to \R$ is some nonlinear, non-convex objective function and $\EE_\PP$ indicates the mathematical expectation with respect to $\PP$, that is, for any $x \in \Rd$,
\[
\EE_\PP[F(x,\mathbf{Y})] = \int_{\Omega} F(x,\mathbf{Y}(\omega)) d\PP(\omega) = \int_E F(x,\mathbf{y}) d\nuy(\mathbf{y}).
\]
In the manuscript, we will assume that $\mathbf{y} \mapsto F(x,\mathbf{y})$ is measurable for all $x$ and positive, $\EE_\PP[F(x,\mathbf{Y})]$ is finite for all $x \in \Rd$, i.e.,
\[ F(x,\mathbf{Y}(\cdot)) \in L^1(\Omega, \mathcal{A}, \PP) \quad \tn{or equivalently} \quad F(x,\cdot) \in L^1_{\nu^{\mathbf{Y}}}(E, \B(E)).
\]
In addition, we require that $f:\Rd \to \R$ admits a global minimizer $\xs \in \Rd$. \\

Among the approaches for solving \eqref{eqi: min problem main}, we are interested in meta-heuristics in order to treat possibly non-differentiable and non-convex functions $F$. These are exploring the search space by balancing the exploitation of the accumulated search experience and exploration \cite{blum2003metaheuristics, bianchi2009survey}. 
Some notable examples include Ant Colony Optimization, Genetic Algorithms, Particle Swarm Optimization, Consensus-Based Optimization and Simulated Annealing.
The growth of interest in such procedures has been a result of their ability of finding a nearly optimal solution for problem instances of realistic size in a generally reasonable computation time, see e.g. \cite{gutjahr2000stochastic, gutjahr2004s,bianchi2004metaheuristics}.
Despite this benefit, most of these methods lack a formal mathematical justification for their efficiency. 

Recent progress in this direction has been made for the class of Consensus-Based Optimization (CBO) algorithms being amenable to theoretical analysis via tools from statistical physics \cite{pinnau2017consensus,carrillo2018analytical,carrillo2021consensus,fornasier2021consensus,fornasier2022convergence}.
These methods use a set of $N \in \mathbb{N}$ particles to explore the domain and find a global minimum of the objective function $f$. The dynamics of the agents is governed by a drift term that drags them towards a momentaneous consensus point (that serves as a temporary guess for $\xs$) and a diffusion term, that favors the exploration of the search space. 
Due to the derivative-free nature of the CBO algorithms, they are applied to (non-convex) non-smooth functions; furthermore, several variations of such methods have been proposed to accommodate for the presence of constraints \cite{borghi2023constrained,carrillo2023consensus,fornasier2022anisotropic,fornasier2021consensusspehere}, multiple objectives \cite{borghi2022consensus, borghi2023adaptive}, objectives for which only a stochastic estimator at a given point can be computed \cite{bellavia2024noisy}
and multiple minimizers \cite{fornasier2024pde}.\\
The theoretical analysis of the CBO methods may be conducted at the algorithmic (hereafter called microscopic) level, as done for instance in \cite{ha2020convergence, ha2021convergence, ko2022convergence,bellavia2024noisy}, or at the mean-field level \cite{bolley2011stochastic} where a statistical description of the dynamics of the average agent behavior is considered. 
The latter viewpoint was proposed in \cite{pinnau2017consensus,carrillo2018analytical} and further applied in many other \rev{works} (\cite{carrillo2021consensus,fornasier2021consensus,fornasier2022convergence,fornasier2022anisotropic} to name a few). In the following we adopt the latter perspective.\\

In this manuscript, we develop two approaches for tackling problem \eqref{eqi: min problem main}. Common to both methods is the approximation of the true objective function $f$ by a suitable sampling and the resolution of the newly obtained problem with the designated CBO algorithm.
The first approach is based on the Sample Average Approximation (SAA) \cite{shapiro2003monte, shapiro2021lectures} and it consists of substituting $f$ with a Monte Carlo estimator $\hat{f}_M$ that depends on a sample $(\mathbf{y}^{(1)}, \ldots, \mathbf{y}^{(M)}) \in E^M$ of $M$ realizations of the random vector $\mathbf{Y}$.
Similarly to the previous works on CBO, we easily derive the mean-field equations for the microscopic systems associated to the objective functions $f$ and $\hat{f}_M$ in the limit $N \to \infty$ (horizontal arrows of Figure \ref{diagi: 2 diagrams}). However, it is not clear a priori, what is the connection between the two just-mentioned mean-field formulations: in this paper we investigate this question and prove that two relations can be established between them (right vertical arrow of Figure \ref{diagi: 2 diagrams}).
The second approach (hereafter called quadrature approach) consists of approximating $f$ by a suitable quadrature formula $\tilde{f}_N$ with $M=N$ fixed discretization points $(\mathbf{y}^1, \ldots, \mathbf{y}^N) \in E^N$.
We require the number of agents in the algorithm to be equal to the number of nodes $N$: this allows us to consider the dynamics of the agents to take place in the augmented space $\Rd \times E$ and to derive the corresponding mean-field equation on the extended phase space.
We then verify that this approach yields the same limiting equation as the one obtained for the SAA (inner diagonal of Figure \ref{diagi: 2 diagrams}). 

\tikzstyle{CBOalg} = [rectangle, rounded corners, minimum width=6cm, minimum height=1cm,text centered, text width=6cm, draw=black, fill=pink!30]
\tikzstyle{MF} = [rectangle, rounded corners, minimum width=7cm, minimum height=1cm,text centered, text width=7cm, draw=black]
\tikzstyle{arrow} = [thick,->,>=stealth]
\tikzstyle{MF2} = [rectangle, rounded corners, minimum width=8cm, minimum height=1cm,text centered, text width=8cm, draw=black, fill=yellow!20]
\begin{figure}[h!tbp]
	\centering
	\resizebox{13cm}{2.7cm}{%
		\begin{tikzpicture}[node distance=2cm]
			\node (cbof) [CBOalg] {CBO algorithm for $f$};
			\node (cboM) [CBOalg, below of=cbof] {CBO algorithm for\\ discretized $f$ (either $\hat{f}_M$ or $\tilde{f}_N$)};
			\node (MFf) [MF, right of=cbof, xshift=6.5cm] {Mean-field formulation for $f$};
			\node (MFM) [MF, right of=cboM, xshift=6.5cm] {Mean-field formulation for \\ discretized $f$};
			
			\draw [arrow,dashed,very thick] (cbof) -- node[anchor=north]{$N \to \infty$}(MFf);
			\draw [arrow,dashed,very thick] (cboM) -- node[anchor=north]{$N \to \infty$} (MFM);
			\draw [arrow] (cbof) -- node[anchor=west]{discretization} (cboM);
			\draw [arrow,very thick,<->] (MFM) -- node[anchor=west]{$M \to \infty$ (*)} (MFf);
			\draw [arrow,blue,very thick,dashdotted] (cboM) -- node[anchor=west,xshift=0.5cm,yshift=-0.2cm]{$N=M \to \infty$}(MFf);
	\end{tikzpicture}}
	\caption{Diagram illustrating the limits derived in the manuscript for the SAA (outer loop) and quadrature (inner diagonal) approaches. The precise sense in which the convergences apply will be specified in Sections \ref{sec:SAA} and \ref{sec:quadr}. The asterisk (*) denotes the fact that additional variables are sent to infinity to obtain the vertical arrow in the diagram.}
	\label{diagi: 2 diagrams}
\end{figure}

\newpage
The rest of the paper is organized as follows. 
In Section \ref{sec:SAA}, we illustrate the combination of the SAA approach and of the CBO algorithm. After presenting the algorithm, we give a detailed proof of the arrows of the outer loop of Figure \ref{diagi: 2 diagrams}. We justify the hypotheses made in view of the key theoretical properties of the solutions obtained through the SAA approach. In Section \ref{sec:quadr}, we present the quadrature approach and verify the equality of the mean-field formulations of the two approaches, hence, obtaining the diagonal arrow of Figure \ref{diagi: 2 diagrams}. We devote Section \ref{sec:numerics} to some numerical tests aiming at investigating the rates of the SAA and quadrature approach and at comparing the two methods. We summarize our main conclusions in Section \ref{sec: saa+cbo vs quad+cbo} and provide an overview of possible directions for further research.

\section{The sample average approximation approach}
\label{sec:SAA}

A common approach for tackling global optimization problem \eqref{eqi: min problem main} is the Sample Average Approximation (SAA) approach (see e.g. \cite{shapiro2003monte, shapiro2021lectures}). It consists of three steps:  
\begin{enumerate}
	\item sample $\bY(\cdot) := (\mathbf{Y}^{(1)}(\cdot), \ldots, \mathbf{Y}^{(M)}(\cdot))$ from the random vector $\mathbf{Y}$, for a fixed $M \in \mathbb{N}$;
	\item approximate $f(x) = \EE_\PP[F(x,\mathbf{Y})]$ with the so-called Sample Average Approximation (SAA) 
	\begin{equation}
		\label{def2: fhatM}
		\hat{f}_M(x, \bY(\cdot)) := \frac{1}{M} \sum_{j=1}^M F(x, \mathbf{Y}^{(j)}(\cdot)) \quad \tn{for} \; x \in \Rd.
	\end{equation}
	\item solve the optimization problem
	\begin{equation}
		\label{eq2: min problem, saa}
		\min_{x \in \Rd} \hat{f}_M(x, \bY(\cdot)).
	\end{equation}
\end{enumerate}
\rev{In the definition of the SAA approach, there is no requirement for independence among the set $\{\bfY^{(j)}(\cdot)\}^j$: we will elaborate on this point in Remark \ref{rem:unifassumpSAA}.}
As for the true problem \eqref{eqi: min problem main}, we require that \eqref{eq2: min problem, saa} admits a unique minimizer $\hat{x}_M(\bY(\cdot))$.
\rev{In the following, we use the notation
	\begin{equation}
		\label{def:fbar,fMbar}
		\underline{f} := \inf_{x} f(x) = f(\xs), \quad \underline{\hat{f}_M}(\bY(\cdot)) := \inf_{x} \hat{f}_M(x,\bY(\cdot)) = \hat{f}_M(\hat{x}_M(\bY(\cdot)),\bY(\cdot)),
	\end{equation}
	where the two equalities are guaranteed by the uniqueness of the minimizers $\xs$ and $\hat{x}_M(\bY(\cdot))$ respectively.}
The notations $\bY(\cdot), \bfY^{(j)}(\cdot)$ aim at emphasizing that $\{\hat{f}_M(x, \bY(\cdot))\}_M$, for a given $x \in \Rd$, is a sequence of random functions defined on the common probability space $(\Omega, \A, \PP)$. 
As problem \eqref{eq2: min problem, saa} is a function of the considered sample $\bY(\cdot)$, it is in that sense random \cite{shapiro2021lectures}.
Subsequently, for a given realization $\omega \in \Omega$, we introduce the notation
$\by = (\mathbf{y}^{(1)}, \ldots,\mathbf{y}^{(M)}) := (\mathbf{Y}^{(1)}(\omega), \ldots, \mathbf{Y}^{(M)}(\omega)) \in E^M,$
that is, we use lower case letters for the $M$ realizations of the random vector $\bfY$.
In particular, $x \mapsto \hat{f}_M(x, \bY(\omega)) = \hat{f}_M(x, \by)$ becomes a deterministic function and \eqref{eq2: min problem, saa} reduces to a deterministic optimization problem.\\

Thanks to the just-mentioned approach, we then have two optimization problems: \eqref{eqi: min problem main} for the stochastic objective $f$ and \eqref{eq2: min problem, saa} for its approximation $\hat{f}_M(\bY(\cdot))$.
In the following, we fix a time horizon $T>0$ and apply two CBO-type algorithms to solve problems \eqref{eqi: min problem main} and \eqref{eq2: min problem, saa} respectively. 

For problem \eqref{eqi: min problem main}, we consider a stochastic system of $N \in \mathbb{N}$ agents with position vectors $X^i_t \in \Rd$, $i=1,\ldots,N$, that dynamically interact with each other to find the minimizer of the objective function.
More precisely, the CBO update rule is given, at time $t \in [0,T]$, by the system of SDEs
\begin{subequations}
	\label{eq2: complete cbo for f, limMinfty done}
	\begin{align}
		dX^i_t &= -\lambda(X^i_t - x^{\alpha,f}_t) dt + \sigma D^i_t dB^i_t \quad \textnormal{for} \; i=1,\ldots,N, \label{eq2: sde cbo for EF}\\
		x^{\alpha,f}_t &= \quad \frac{\sum_{i=1}^{N} X^i_t \exp(-\alpha f(X^i_t))}{\sum_{i=1}^{N} \exp(-\alpha f(X^i_t))},  \label{eq2: xalphat_f}
	\end{align}
\end{subequations}
where $\lambda, \sigma >0$, $\alpha >0$, and $\{B^i_t\}^i$ are $d$-dimensional independent Brownian processes. 
In system \eqref{eq2: sde cbo for EF} the drift term is governed by $\lambda$, moving particles towards the consensus point $x^{\alpha,f}_t$, and the diffusion term is governed by $\sigma$, for exploration in the search space. $D^i_t$ is a matrix in $\R^{d \times d}$ that can be chosen either as 
\[
D^i_t = D^i_{t,\textnormal{iso}} := |X^i_t - x^{\alpha,f}_t|I_d, \quad \tn{where $I_d \in \R^{d \times d}$ is the identity matrix},
\]
so that the random exploration process is isotropic \cite{pinnau2017consensus} and all dimensions $l = 1, \ldots, d$ are equally explored, or as 
\[
D^i_t = D^i_{t,\textnormal{aniso}} :=  \tn{diag}(X^i_t -x^{\alpha,f}_t),
\]
so that it is anisotropic \cite{carrillo2021consensus}. 
The exponential coefficients in the weighted average $x^{\alpha,f}_t$ \eqref{eq2: xalphat_f} fulfill
\[
x^{\alpha,f}_t \xrightarrow[]{\alpha \to \infty} \tn{argmin}_{i=1,\ldots,N} f(X^i_t).
\]
The above convergence is guaranteed by the Laplace principle \cite{dembo2009large}, see also \cite{pinnau2017consensus,carrillo2018analytical,carrillo2021consensus,fornasier2021consensus,fornasier2022convergence}. 
The system is supplemented with initial conditions $X^i_0, i=1,\ldots,N$, independent and identically distributed (i.i.d.) with law $\mu_0 \in \mathcal{P}(\Rd)$. 

Subsequently, for problem \eqref{eq2: min problem, saa}, 
we consider a stochastic system of $N \in \mathbb{N}$ agents with position vectors $X^i_t(\bY(\cdot)) \in \Rd$ and describe the update rule by
\begin{subequations}
	\label{eq2: complete cbo for fhatM}
	\begin{align}
		dX^i_t(\bY(\cdot)) &= -\lambda(X^i_t(\bY(\cdot)) - x^{\alpha,\hat{f}_M(\bY(\cdot))}_t) dt + \sigma D^i_t(\bY(\cdot)) dB^i_t \quad \textnormal{for} \; i, \label{eq2: sde cbo for fhatM}\\
		x^{\alpha,\hat{f}_M(\bY(\cdot))}_t &= \frac{\sum_{i=1}^{N} X^i_t(\bY(\cdot)) \exp(-\alpha \hat{f}_M(X^i_t(\bY(\cdot)), \bY(\cdot)))}{\sum_{i=1}^{N} \exp(-\alpha \hat{f}_M(X^i_t(\bY(\cdot)), \bY(\cdot)))}, \label{eq2: xalphat_M}
	\end{align}
\end{subequations}
with the same parameters $\lambda, \sigma, \alpha$ and initial conditions $\{X^i_0\}^i$ as above, and $D^i_t(\bY(\cdot))$ equal to 
\begin{equation*}
	D^i_{t,\tn{iso}}(\bY) = |X^i_t(\bY(\cdot)) - x^{\alpha,\hat{f}_M(\bY(\cdot))}_t| I_d,  D^i_{t,\tn{aniso}}(\bY) = \tn{diag}(X^i_t(\bY(\cdot)) - x^{\alpha,\hat{f}_M(\bY(\cdot))}_t).
\end{equation*}
\noindent For the sake of simplicity, we limit our theoretical analysis to the case where particle positions are updated using isotropic diffusion (i.e., $D^i_t = D^i_{t,\tn{iso}},D^i_{t,\tn{iso}}(\bY(\cdot))$).\\

It is straightforward to verify that 
\begin{assumption}
	\label{ass:well-posmicro_F}
	$F(x,\bfy)$ is locally Lipschitz continuous in the first variable $x \in \Rd$, uniformly in the second variable $\bfy \in E$. 
	More precisely, for all $n >0$, there exists $L_F > 0$ (possibly dependent on $n$) such that
	\begin{equation*}
		|F(x_1,\bfy) - F(x_2,\bfy) | \le L_F |x_1-x_2| \quad \tn{$\forall \; (x_1,x_2,\bfy) \in \R^{2d} \times E$, with $|x_1|,|x_2| \le n$.}
	\end{equation*}
\end{assumption}
implies that $f$ and $x \mapsto \hat{f}_M(x, \bY(\omega))$, for a given realization $\omega \in \Omega$, are locally Lipschitz continuous with the same Lipschitz constant $L_F$. Hence, provided that also $(X^1_0,\ldots,X^N_0)$ has finite second moment, and observing that $f$ and $\hat{f}_M$ are lower bounded by $\underline{f}, \underline{\hat{f}}_M$ defined in \eqref{def:fbar,fMbar}, it holds that
\eqref{eq2: complete cbo for f, limMinfty done} and \eqref{eq2: complete cbo for fhatM} admit unique strong solutions $\{X^i_t\}^i_t$ and $\{X^i_t(\bY(\omega))\}^i_t$ respectively \cite[Theorem 2.1]{carrillo2018analytical}, \cite[Theorem 3.1]{durrett1996stochastic}, defined on suitable probability spaces (see also Remark \ref{rem:probspace}).

\begin{remark}
	\label{rem:probspace}
	Both $X^i_t$ and $X^i_t(\bY(\omega))$, for a fixed $\omega \in \Omega$, are solutions to SDEs, hence are random vectors defined on appropriate probability spaces. 
	We denote the spaces by $(\tilde{\Omega}, \tilde{\mathcal{A}}, \tilde{\PP})$ and $(\dbtilde{\Omega}, \dbtilde{\mathcal{A}}, \dbtilde{\PP})$ respectively, and observe that they are distinct.
	The notation $X^i_t, X^i_t(\bY(\cdot)) \in \Rd$ should then be intended as $X^i_t(\tilde{\omega}) \in \Rd$ for any $\tilde{\omega} \in \tilde{\Omega}$ and $X^i_t(\bY(\omega))(\dbtilde{\omega}) \in \Rd$ for any $\dbtilde{\omega} \in \dbtilde{\Omega}$. In particular, the expression $X^i_t(\bY(\omega))(\dbtilde{\omega})$ aims at emphasizing that the randomness of the random vector $\bY$ is separated from that of the algorithm.
\end{remark}

\subsection{The mean-field equations for $N \to \infty$}
\label{subsec:SAA_MFlevelder}

Similarly to the previous works on CBO, we formallly derive the mean-field equation associated to \eqref{eq2: complete cbo for f, limMinfty done} in the limit $N \to \infty$ by using the propagation of chaos assumption on the marginals \cite{cercignani2013mathematical,golse2003mean,jabin2017mean}.
For a fixed time $t$, let $H^N(t)$ be the $N$-particle probability distribution associated to \eqref{eq2: complete cbo for f, limMinfty done}. Then, the assumption translates to  \rev{$H^N(t) \approx h(t)^{\otimes N}$} for a probability distribution $h(t)$ over $\Rd$.
Finally, we obtain that $h$ \rev{formally} satisfies in the distributional sense
\begin{subequations}
	\label{eq2: mf eq complete for EF}
	\begin{align}
		\pd_t h(t,x) &= \lambda \nabla_x \cdot \left(  (x - x^{\alpha, f}_t[h] ) \:h(t,x) \right) + \frac{\sigma^2}{2} \Delta_x \left( |x - x^{\alpha, f}_t[h]|^2 \:h(t,x) \right),
		\label{eq2: mf eq for EF}\\
		x^{\alpha, f}_t[h] &= \frac{\int_{\Rd} x \exp(-\alpha f(x)) h(t,x) dx}{\int_{\Rd} \exp(-\alpha f(x)) h(t,x) dx}, \label{eq2: xalphat_f_bar}\\
		\mu_0(dx) &=
		\lim_{t \to 0} h(t,x) dx.  
	\end{align}
\end{subequations}

Letting $H^N(t)(\bY(\cdot))$ be the $N$-particle probability distribution at time $t$ associated to \eqref{eq2: complete cbo for fhatM}, we apply again the propagation of chaos assumption to \rev{formally} conclude that \rev{$H^N(t)(\bY(\cdot)) \approx h(t)^{\otimes N}(\bY(\cdot))$}, for a random probability distribution $h(t)(\bY(\cdot))$ over $\Rd$, and $h(\bY(\cdot))$ satisfies in the distributional sense 
\begin{subequations}
	\label{eq2: mf eq complete for fM}
	\begin{align}
		\pd_t h(t,x)(\bY(\cdot)) &= 
		\lambda \nabla_x \cdot \left(  (x - x^{\alpha, \hat{f}_M(\bY(\cdot))}_t[h(\bY(\cdot))] ) \:h(t,x)(\bY(\cdot)) \right) \nonumber\\
		&\hspace{1cm} + \frac{\sigma^2}{2} \Delta_x \left( \left|x-x^{\alpha, \hat{f}_M(\bY(\cdot))}_t[h(\bY(\cdot))] \right|^2 \:h(t,x) \right),
		\label{eq2: mf eq for fM}\\
		x^{\alpha, \hat{f}_M(\bY(\cdot))}_t[h(\bY(\cdot))] &= \frac{\int_{\Rd} x \exp(-\alpha \hat{f}_M(x,\bY(\cdot))) h(t,x)(\bY(\cdot)) dx}{\int_{\Rd} \exp(-\alpha \hat{f}_M(x,\bY(\cdot))) h(t,x)(\bY(\cdot)) dx}, \label{eq2: xalphat_M_bar}\\
		\mu_0(dx) &=
		\lim_{t \to 0} h(t,x)(\bY(\cdot)) dx.
	\end{align}
\end{subequations}
\noindent We remark that the presence of the random quantity $\bY(\cdot)$ does not influence the passage to the mean-field limit as the derivation can be carried out for a fixed realization $\omega \in \Omega$.\\

We denote by $|\cdot|$ the Euclidean norm in $\Rd$. 
We make boundedness and smoothness assumptions regarding our cost function $F$ (we adapt the assumptions of \rev{\cite[Theorem 2.1]{carrillo2018analytical}, \cite[Theorems 3.1, 3.2]{carrillo2018analytical}} to $F$, function of two variables).
\begin{assumption} 
	\label{ass:well-posMF_F}
	\begin{enumerate}
		\item For any $\bfy \in E$, $x \mapsto F(x,\bfy)$ is bounded from below \rev{by} $\underline{F}(\bfy):= \inf_{x \in \Rd} F(x,\bfy)$.
		\item There exist constants $J_F, c_u >0$ such that, for any $\bfy \in E$,
		\begin{equation*}
			\begin{cases}
				|F(x_1,\bfy) - F(x_2,\bfy)| \le J_F (\rev{1+} |x_1|+|x_2|)|x_1-x_2| \quad \tn{for all $x_1,x_2 \in \Rd$}, \\
				F(x,\bfy) - \underline{F}(\bfy) \le c_u (1+|x|^2) \quad \tn{for all $x \in \Rd$.}
			\end{cases}
		\end{equation*}
		\item \rev{Either one of the two conditions is true:}
		\begin{itemize}
			\item \rev{for} any $\bfy \in E$, $x \mapsto F(x,\bfy)$ is bounded \rev{by} $\overline{F}(\bfy):= \sup_{x \in \Rd} F(x,\bfy)$ from above,
			\item \emph{or} there exist $\bar{c_l}>0$ and, for $\bfy \in E$, there exist $c_l(\bfy)>0$ such that
			\begin{equation*}
				F(x,\bfy) - \underline{F}(\bfy) \ge c_l(\bfy) |x|^2 \quad \tn{for all $|x| \ge \bar{c_l}$.}
			\end{equation*}
		\end{itemize}
	\end{enumerate}
\end{assumption}
Assumption \ref{ass:well-posMF_F} implies that $f$ and $x \mapsto \hat{f}_M(x, \bY(\omega))$, for a given realization $\omega \in \Omega$,  satisfy the \rev{requirements} of \rev{\cite[Theorems 3.1, 3.2]{carrillo2018analytical}}: then, under the additional \rev{condition} $\mu_0 \in \mathcal{P}_4(\Rd)$, the aforementioned theorems guarantee the existence of a weak solution $h(t,x)dx \in \mathcal{P}_4(\Rd)$ of \eqref{eq2: mf eq complete for EF} and $h(t,x)(\bY(\omega))dx \in \mathcal{P}_4(\Rd)$ of \eqref{eq2: mf eq complete for fM} for $t \in [0,T]$.\\


All things considered, so far we have guaranteed the well-posedness of the solutions to \eqref{eq2: complete cbo for f, limMinfty done} and \eqref{eq2: mf eq complete for EF}, and \eqref{eq2: complete cbo for fhatM} and \eqref{eq2: mf eq complete for fM} for a given realization $\omega \in \Omega$, and formally derived the implications of Figure \ref{diag:onlyN} for fixed $M \in \mathbb{N}$ (and $\bY(\cdot)$), $\alpha>0$ and $t \in [0,T]$, under the \rev{assumption} of the validity of the propagation of chaos assumption on the marginals. We note that the above formal derivation may be made rigorous by proceeding as in \cite{huang2022mean} (in particular, for the lower arrow it is sufficient to fix a realization $\omega \in \Omega$ to be in the setting required by \cite{huang2022mean}).
\noindent In Subsection \ref{subsec:SAA_MFlevel}, we investigate the completion of the diagram of Figure \ref{diag:onlyN} with a right vertical arrow, namely we investigate whether we can establish a connection between mean-field formulations \eqref{eq2: mf eq complete for EF} and \eqref{eq2: mf eq complete for fM}.
\begin{figure}[h!tbp]
	\centering
	\resizebox{13cm}{2.7cm}{%
		\begin{tikzpicture}[node distance=2cm]	
			\node (cboM) [CBOalg] {CBO algorithm \eqref{eq2: complete cbo for fhatM}\\ with consensus $x^{\alpha,\hat{f}_M(\bY(\cdot))}_t$ \eqref{eq2: xalphat_M}};
			\node (cbof) [CBOalg, above of=cboM] {CBO algorithm \eqref{eq2: complete cbo for f, limMinfty done}\\ with consensus $x^{\alpha, f}_t$ \eqref{eq2: xalphat_f}};
			\node (MFM) [MF, right of=cboM, xshift=6.5cm] {Mean-field formulation \eqref{eq2: mf eq complete for fM}\\with consensus $x^{\alpha, \hat{f}_M(\bY(\cdot))}_t[h(\bY(\cdot))]$ \eqref{eq2: xalphat_M_bar}};
			\node (MFf) [MF, above of=MFM] {Mean-field formulation \eqref{eq2: mf eq complete for EF}\\with consensus $x^{\alpha, f}_t[h]$ \eqref{eq2: xalphat_f_bar}};
			\draw [arrow,dashed,very thick] (cboM) -- node[anchor=north]{$N \to \infty$} (MFM);
			\draw [arrow,dashed,very thick] (cbof) -- node[anchor=north]{$N \to \infty$}(MFf);
	\end{tikzpicture}}
	\caption{Derivation of mean-field formulations  \eqref{eq2: mf eq complete for EF} and \eqref{eq2: mf eq complete for fM} in the limit of the number of agents $N$ going to infinity and for fixed $M \in \mathbb{N}$ (and $\bY(\cdot)$), $\alpha>0$ and $t \in [0,T]$. The sense in which the convergence holds is the one that can be formally derived through the propagation of chaos assumption on the marginals and that can be made rigorous by proceeding as in \cite{huang2022mean}. For the bottom arrow, the convergence holds for any fixed realization $\omega \in \Omega$.}
	\label{diag:onlyN}
\end{figure}

\subsection{The limit for $M \to \infty$ at the mean-field level}
\label{subsec:SAA_MFlevel}

In this subsection, we establish two relations between mean-field formulations \eqref{eq2: mf eq complete for EF} and \eqref{eq2: mf eq complete for fM}. We first prove that the $2$-Wasserstein distance  between the solutions of the formulations converges to zero for $M,\alpha,t$ sufficiently big and with probability one (w.p.1) (see Proposition \ref{prop:convsol_MF}; for an introduction on $p$-Wasserstein distances, see e.g. \cite{ambrosio2005gradient,chaintron2022propagation}); then, in the same regime and w.p.1, we prove the stronger result of convergence to zero of the Euclidean norm of the difference between the corresponding consensus points (see Theorem \ref{th:convcons_MF}). \\

In the following, we denote by $\C(A,B)$ the space of continuous functions from $A$ to $B$.
We make some tractability conditions of the landscape of $f$ and $\hat{f}_M$ around the unique minimizers $\xs, \hat{x}_M$ resp. and in the farfield (we adapt the assumptions of \rev{\cite[Theorem 3.7]{fornasier2021consensus}} to our objectives) and derive our first result in the $2$-Wasserstein distance.
\begin{assumption}
	\label{ass:tract}
	\begin{enumerate}
		\item There exist $f_{\infty},R_0,\eta>0$ and $\nu \in (0,\infty)$ such that
		\begin{equation*}
			\begin{cases}
				|x-\xs| \le \frac{1}{\eta} (f(x)-\underline{f})^{\nu} \quad \tn{for all $|x-\xs| < R_0$} \\
				f(x) - \underline{f} > f_{\infty} \quad \tn{for all $|x-\xs| > R_0$.}
			\end{cases}
		\end{equation*}
		\item For any $\by \in E^M$, there exist $\tilde{f}_{\infty},\tilde{R}_0,\tilde{\eta}>0$ and $\tilde{\nu} \in (0,\infty)$ (all \rev{possibly} depending on $\by$) such that 
		\begin{equation*}
			\begin{cases}
				|x-\hat{x}_M(\by)| \le \frac{1}{\tilde{\eta}} (\hat{f}_M(x,\by)-\underline{\hat{f}_M}(\by))^{\tilde{\nu}} \quad \tn{for all $|x-\hat{x}_M(\by)| < \tilde{R}_0$} \\
				\hat{f}_M(x,\by) - \underline{\hat{f}_M}(\by) > \tilde{f}_{\infty} \quad \tn{for all $|x-\hat{x}_M(\by)| > \tilde{R}_0$.}
			\end{cases}
		\end{equation*}
		for $\underline{f}, \underline{\hat{f}_M}(\bY(\cdot))$ defined in \eqref{def:fbar,fMbar}.
	\end{enumerate}
\end{assumption}

\begin{proposition} 
	\label{prop:convsol_MF}
	Let $\mu_0 \in \mathcal{P}_4(\Rd)$ and let $F(\cdot,\bfy) \in \C(\Rd,\R)$, for any $\bfy \in E$, satisfy Assumptions \ref{ass:well-posMF_F}. 
	Let  $f, \hat{f}_M$ fulfill Assumption \ref{ass:tract} and choose $\lambda,\sigma>0$ with $2\lambda>d\sigma^2$.
	Assume that there exists a nonempty compact subset $C$ of $\Rd$ such that $\xs, \hat{x}_M(\bY(\cdot)) \in \tn{supp}(\mu_0) \subset C$ and such that $\hat{f}_M$ converges to $f$ w.p.1 uniformly on $C$ as $M$ tends to infinity.
	Then, denoting by $h(t,x)dx$ and $h(t,x)(\bY(\cdot))dx$ the weak solutions to mean-field formulations \eqref{eq2: mf eq complete for EF} and \eqref{eq2: mf eq complete for fM} resp., it holds that 
	\begin{equation*}
		W_2(h(t,x)dx,h(t,x)(\bY(\cdot))dx) \xrightarrow[]{} 0 \quad \tn{for $\alpha,t,M$ sufficiently large and w.p.1.}
	\end{equation*}
\end{proposition}

\begin{proof}[Proof of Proposition \ref{prop:convsol_MF}]
	Let $C$ be the compact subset satisfying the assumptions of the theorem. 
	The \rev{requirements} $F(\cdot,\bfy) \in \C(\Rd)$, for any $\bfy \in E$, $\xs, \hat{x}_M(\bY(\cdot)) \in C$ and the convergence of $\hat{f}_M$ to $f$ uniformly on $C$ imply that $ \hat{x}_M(\bY(\cdot))$ approaches $\xs$ for $M$ sufficiently big and w.p.1, see \rev{\cite[Theorem 5.3]{shapiro2021lectures}}.
	
	Let $\N \in \A$ be the set of measure zero (with respect to probability measure $\PP$) that guarantees the uniform convergence of $\hat{f}_M$ to $f$ (for more details, see Remark \ref{rem:unifassumpSAA}); let us fix $\omega^u \in \Omega \setminus \N$ (we use the superscript $u$ to denote the fact that there is uniform convergence in correspondence of it) and let us define $\by^u := \bY(\omega^u)$. 
	As already remarked in the introduction to the section, \rev{condition} $\mu_0 \in \mathcal{P}_4(\Rd)$ and Assumption \ref{ass:well-posMF_F} guarantee the existence of a weak solution  $h(t,x)dx$ of \eqref{eq2: mf eq complete for EF} and $h(t,x)(\by^u)dx$ of \eqref{eq2: mf eq complete for fM} for $t \in [0,T]$.
	Let us consider the $2$-Wasserstein distance between the aforementioned solutions, namely $W_2(h(t,x)dx,h(t,x)(\by^u)dx)$. Then, by applying twice a triangular inequality with the measures $\delta_{\xs}$ and $\delta_{\hat{x}_M(\by^u)}$ (with $\delta_a$ denoting the Dirac measure of $a \in \Rd$), we get that 
	\rev{
		\begin{align}
			\label{eq:proofprop_W2triang}
			W_2&(h(t,x)dx,h(t,x)(\by^u)dx)  \\
			&\le \underbrace{W_2(h(t,x)dx,\delta_{\xs})}_{(i)} + \underbrace{W_2(\delta_{\xs},\delta_{\hat{x}_M(\by^u)} )}_{(ii)} + \underbrace{W_2(\delta_{\hat{x}_M(\by^u)},h(t,x)(\by^u)dx)}_{(iii)}. \nonumber
	\end{align}}
	
	We now study each term separately, starting out from $(ii)$. Using the definition of $2$-Wasserstein distance between two Dirac measures and the fact that $\omega^u$ is chosen in a set of probability 1, it holds that 
	\begin{equation*}
		\label{eq:hatxMtoxs_yu}
		(ii) = W_2(\delta_{\xs},\delta_{\hat{x}_M(\by^u)} ) = |\xs-\hat{x}_M(\by^u)| \xrightarrow[]{M \to \infty} 0.
	\end{equation*}
	The \rev{assumptions} $\xs \in \tn{supp}(\mu_0)$, $f \in \C(\Rd)$ (implied by the corresponding assumption on $F$), on $\lambda,\sigma$, and condition 1. of Assumption \ref{ass:tract}, along with the ones required to guarantee the well-posedness of mean-field measure $h(t,x)dx$, guarantee
	\begin{equation*}
		\label{eq:W2h_yu}
		(iii) = W_2(h(t,x)dx, \delta_{\xs}) \xrightarrow[]{} 0 \quad {\tn{for $\alpha,t$ sufficiently big}},
	\end{equation*}
	see \rev{\cite[Theorem 3.7]{fornasier2021consensus}}.\\
	Analogously, $\hat{x}_M(\by^u) \in \tn{supp}(\mu_0)$, $\hat{f}_M(\cdot,\by^u) \in \C(\Rd)$, $2\lambda > d\sigma^2$ and condition 2. of Assumption \ref{ass:tract} for $\by = \by^u$ imply
	\begin{equation*}
		\label{eq:W2hY_yu}
		(i) = W_2(h(t,x)(\by^u)dx,\delta_{\hat{x}_M(\by^u)}) \xrightarrow[]{} 0 \quad {\tn{for $\alpha,t$ sufficiently big}}.
	\end{equation*}
	
	Finally, putting together $(i),(ii)$ and $(iii)$ leads us to the desired conclusion.
\end{proof}

\begin{remark}
	\label{rem:unifassumpSAA}
	The proof of Proposition \ref{prop:convsol_MF} is based on one of the building blocks of the SAA theory, namely the consistency of the estimator $\hat{x}_M$ of $\xs$. We obtained the consistency by using the result \rev{\cite[Theorem 5.3]{shapiro2021lectures}} and, in particular, by making the assumption that there exists some nonempty compact subset $C \in \Rd$ such that $\hat{f}_M$ converges to $f$ uniformly on $C$. In more detail, the condition translates to the existence of $\N \in \A$ a set of measure zero (with respect to probability measure $\PP$) such that for any $\omega \in \Omega \setminus \N$
	\begin{equation*}
		\lim_{M \to \infty} \sup_{x \in C} | \hat{f}_M(x, \bY(\omega)) - f(x) | = 0.
	\end{equation*}
	Sufficient conditions for the limit to hold are given for instance in \rev{\cite[Theorems 7.48,7.50]{shapiro2021lectures}} and are based on the \rev{requirement} that the sample $\left\lbrace \mathbf{Y}^{(j)} \right\rbrace^{j}$ is i.i.d. with law $\nuy$.
	Kim and colleagues \cite{kim2015guide} remark that the i.i.d. assumption is not essential and it is sometimes relaxed in the context of variance reduction \cite{caflisch1998monte} and quasi-Monte Carlo schemes \cite{niederreiter1992quasi}, where dependence is inserted on purpose within the sample with the aim of reducing the variance of $\hat{f}_M$ or accelerating the rate $O(M^{-1/2})$ for low dimensional random spaces ($k =1,2$) respectively, hence enhancing the convergence of the corresponding SAA estimators $\hat{x}_M, \hat{f}_M(\hat{x}_M )$. 
	In contrast, Wang and colleagues \cite{wang2022sample} affirm that the assumption of i.i.d. samples is invalid in practice in many data-generating processes and state that assuming that the sample is generated by an ergodic stochastic process that converges to the stationary distribution $\nuy$ is more realistic. Even so, the analysis as well as the final results are similar to the ones of the classical theory of \cite{shapiro2021lectures}, hence, in order to gain an initial understanding of uniform convergence theory, it is sufficient to think of the sample as i.i.d..\\
	Last, the consistency of $\hat{x}_M$ can be obtained by replacing the uniform convergence property with the notion of epiconvergence, see for instance \rev{\cite[Theorem 5.4]{shapiro2021lectures}}. However, this approach restricts to convex objectives $F$ and, as the class of CBO algorithms treats possibly non-differentiable and non-convex function, we decide to avoid it in our analysis.\\
\end{remark}

Now we prove a result stronger than Proposition \ref{prop:convsol_MF} that states that the Euclidean norm of the difference between consensus points \eqref{eq2: xalphat_f_bar} and \eqref{eq2: xalphat_M_bar}
converges to zero for $\alpha,t,M$ sufficiently big and w.p.1 (Theorem \ref{th:convcons_MF}). 
We require the randomness of the function $x \mapsto F(x, \mathbf{Y})$ to be confined to a compact subset. An example of function satisfying such assumption is given in Subsection \ref{subsecnum: rate SAA}.
\begin{assumption}
	\label{ass:confrand_F}
	For $x$ sufficiently large, $x \mapsto F(x,\bfy)$ is independent of $\bfy$. In particular, there exist a non empty sufficiently big compact subset $C$ of $\Rd$ such that, for any $x \in \Rd \setminus C$, $x \mapsto F(x,\bfy)$ is independent of $\bfy$.
\end{assumption}
\noindent We premise three lemmas that we will use in the proof of the theorem. 

\begin{lemma} [\rev{\cite[Theorem 7.1.2]{arnold1974stochastic}} for $\mu_0 \in \mathcal{P}_4(\Rd)$]
	\label{lem:4thestimate}
	Suppose that we have a stochastic differential equation
	\begin{equation}
		\label{eq:genSDE}
		dW_t = \Delta(t,W_t) dt + \Sigma(t,W_t) dB_t, \; \tn{for $t \in [0,T]$,} \quad \tn{and} \quad W_0 \sim \mu_0,
	\end{equation}
	for some $\Delta(\cdot,\cdot):[0,T]\times \Rd \to \Rd, \Sigma(\cdot,\cdot):[0,T]\times \Rd \to \Rd$ measurable vector-values and matrix-valued functions respectively, $B_t$ the standard Brownian motion in $\Rd$ and for some fixed initial measure $\mu_0 \in \mathcal{P}_4(\Rd)$. 
	Let $L > 0$ be the Lipschitz constant such that $\Delta$ and $\Sigma$ are locally Lipschitz continuous in the second variable uniformly in the first
	and let $\mu \in \C([0,T], \mathcal{P}_4(\Rd))$ be such that $W_t \sim \mu_t$. 
	Then, the following fourth-order moment estimate for $\mu_t$ holds
	\begin{equation*}
		\int_{\Rd} |x|^4 \mu_t(dx) \le \left(1 + \int_{\Rd} |x|^4 \mu_0(dx)\right) e^{20 L^2 t}.
	\end{equation*}
\end{lemma}

\begin{lemma}[\rev{\cite[Lemma 3.2]{carrillo2018analytical}}]
	\label{lem:estimatecons_diffmeasCCTT18}
	Let $\mathcal{F}: \Rd \to \R$ be an objective function bounded from below with $\underline{\F} := \inf_{x \in \Rd} \F(x)$ and such that there exist constants $J_{\F}, c_{u,\F}>0$ such that
	\begin{equation*}
		\begin{cases}
			|\F(x_1) -\F(x_2)| \le J_{\F} (\rev{1+} |x_1|+|x_2|)|x_1-x_2| \quad \tn{for all $x_1,x_2 \in \Rd$}, \\
			\F(x) - \underline{\F} \le c_{u,\F} (1+|x|^2) \quad \tn{for all $x \in \Rd$.}
		\end{cases}
	\end{equation*}
	Let $\mu, \tilde{\mu} \in \mathcal{P}_2(\Rd)$ with $4$-th moment bounded by a constant $k_4 < \infty$. Then the following stability estimate holds
	\begin{equation*}
		|x^{\alpha,\mathcal{F}}[\mu] - x^{\alpha,\mathcal{F}}[\tilde{\mu}] | \le c_0 W_2(\mu, \tilde{\mu}),
	\end{equation*}
	for a constant $c_0$ depending only on $\alpha >0, J_{\F}$ and $k_4$ and for 
	\[
	x^{\alpha,\mathcal{F}}[\nu] = \frac{\int_{\Rd} x \exp(-\alpha \F(x)) \nu(dx)}{\int_{\Rd} \exp(-\alpha \F(x)) \nu(dx)} \quad \tn{for $\nu = \mu, \tilde{\mu}$}. 
	\]
\end{lemma}

\begin{lemma} 
	\label{lem:convcons_MF}
	Let $\mu_0 \in \mathcal{P}_4(\Rd)$ and let $F$ satisfy Assumptions \ref{ass:well-posmicro_F}, \ref{ass:well-posMF_F} and \ref{ass:confrand_F}. Let $t \in [0,T], M \in \mathbb{N}$ and the parameters $\lambda, \sigma, \alpha > 0$ be fixed. 
	Given $\by \in E^M$ and denoting by $h(t,x)dx$ a weak solution of \eqref{eq2: mf eq complete for EF} and by $h(t,x)(\by)dx$ a weak solution of \eqref{eq2: mf eq complete for fM}, then the following estimates hold: 
	\begin{subequations}
		\begin{align}
			\left| \; x^{\alpha, \hat{f}_M(\by)}_t[h] -  x^{\alpha, \hat{f}_M(\by)}_t[h(\by)] \; \right| &\le C_0 W_2(h(t,x)dx,h(t,x)(\by))dx) \label{eq:estimatexalpha_diffmeas}\\
			\bigg| \int_{\Rd} x \left( e^{-\alpha f(x)} - e^{-\alpha \hat{f}_M(x, \by)} \right)  h(t,x) dx \bigg| &\le 
			C_1 \sup_{x \in C} \left| \; f(x)- \hat{f}_M(x, \by) \; \right| \label{eq:estimatenum_diffobj}\\
			\bigg| \int_{\Rd} \left( e^{-\alpha f(x)} - e^{-\alpha \hat{f}_M(x, \by)} \right)  h(t,x) dx \bigg|  &\le  \sup_{x \in C} \left| \; f(x) - \hat{f}_M(x, \by) \; \right| \label{eq:estimateden_diffobj}
		\end{align} 
	\end{subequations}
	for some constants $C_0$ and $C_1$ both independent of $M$ and $\by$.
\end{lemma}

\begin{proof}[Proof of Lemma \ref{lem:convcons_MF}]
	Let $\by \in E$ be fixed. 
	We begin with the proof of estimate \eqref{eq:estimatexalpha_diffmeas}. As already remarked in the introduction to the section and in Proposition \ref{prop:convsol_MF}, the \rev{condition} $\mu_0 \in \mathcal{P}_4(\Rd)$ and Assumption \ref{ass:well-posMF_F} guarantee the existence of a weak solution $h(t,x)dx \in \mathcal{P}_4(\Rd)$ of \eqref{eq2: mf eq complete for EF} and $h(t,x)(\by)dx \in \mathcal{P}_4(\Rd)$ of \eqref{eq2: mf eq complete for fM} for $t \in [0,T]$.
	As each of the two measures corresponds to a system of SDEs of the form \eqref{eq:genSDE} for suitable choices of $\Delta,\Sigma$, it is possible to apply Lemma \ref{lem:4thestimate} and conclude that 
	\begin{equation}
		\label{eq:4thestimate}
		\int_{\Rd} |x|^4 h(t,x)dx, \int_{\Rd} |x|^4 h(t,x)(\by)dx \le \left(1 + \int_{\Rd} |x|^4 \mu_0(dx)\right) e^{20 L^2 t},
	\end{equation}
	for some constant $L>0$ depending uniquely on $L_F$ defined in Assumption \ref{ass:well-posmicro_F} (hereafter denoted as $L(L_F)$). We set
	\begin{equation}
		\label{def:K4}
		K_4(L_F) :=  \left(1 + \int_{\Rd} |x|^4 \mu_0(dx)\right) e^{20 L(L_F)^2 t}
	\end{equation}
	and remark that it is independent of $\by$ thanks to the independence of $L_F$ of $\by$. 
	Then, we may apply Lemma \ref{lem:estimatecons_diffmeasCCTT18} to $h(t,x)dx, h(t,x)(\by)dx$ and $\hat{f}_M(\cdot,\by)$ by commenting that the latter function satisfies the desired properties with constants $J_F,c_u$ as a consequence of Assumption \ref{ass:well-posMF_F}. In particular, equation \eqref{eq:estimatexalpha_diffmeas} holds for some constant $C_0 >0$ depending uniquely on $\alpha, J_F$ and $K_4(L_F)$ (given by \eqref{def:K4}), all independent of $\by$.
	
	Now we prove \eqref{eq:estimatenum_diffobj}:
	\begin{align*}
		\bigg| &\int_{\Rd} x \left( e^{-\alpha f(x)} - e^{-\alpha \hat{f}_M(x, \by)} \right)  h(t,x) dx \bigg| 
		\le \; \int_{\Rd}\left| x  \left(e^{-\alpha f(x)} - e^{-\alpha \hat{f}_M(x, \by)} \right) \right| h(t,x) dx \\
		&\hspace{0.2cm}\underbrace{\le}_{\rev{\alpha}-\tn{Lip}} \int_{\Rd} |x| \: \rev{\alpha} \left| f(x) - \hat{f}_M(x, \by)  \right| h(t,x) dx\\
		&\hspace{0.4cm}=  \rev{\alpha} \int_{C} |x| \left| f(x) - \hat{f}_M(x, \by) \right| h(t,x) dx + \rev{\alpha} \int_{\Rd \setminus C} |x| \underbrace{\left| f(x) - \hat{f}_M(x, \by) \right|}_{\tn{$=0$ on $\Rd \setminus C$}} h(t,x) dx\\
		&\hspace{0.4cm}= \rev{\alpha} \int_{C} |x| \left| f(x) - \hat{f}_M(x, \by) \right| h(t,x) dx \\
		&\hspace{0.4cm}\le  \sup_{x \in C} | f(x) - \hat{f}_M(x, \by)| \; \rev{\alpha} \: \int_C |x| h(t,x) dx,\\
		&\hspace{-0.3cm}
		\underbrace{\le}_{L^4(C) \subset L^1(C)}  \sup_{x \in C} | f(x) - \hat{f}_M(x, \by)| \; \rev{\alpha} \: \tn{meas}(C)^{\frac{1}{4}} K_4(L_F),
	\end{align*}	
	where in the second row we have used the fact that $ z \mapsto e^{-\alpha z}$ is \rev{$\alpha$}-Lipschitz for $z \ge 0$ (and $\hat{f}_M, f$ \rev{are non-negative due to the corresponding assumption on $F$ made in Section \ref{sec:intro}}), in the subsequent equality that, thanks to Assumption \ref{ass:confrand_F}, for any $x \in \Rd \setminus C$, $\hat{f}_M(x,\by) = f(x)$ and, in the last passage, that $L^4(C) \subset L^1(C)$ (with $L^p$ the Lebesgue space for $p>0$) and that estimate \eqref{eq:4thestimate} holds for $h$. We set $C_1 := \rev{\alpha} \: \tn{meas}(C)^{\frac{1}{4}} K_4(L_F) > 0$ and observe that it is independent of $\by$.
	
	The proof of \eqref{eq:estimateden_diffobj} is analogous to that presented above, except for the last inequality where we instead use that
	\begin{equation*}
		\int_C  h(t,x) dx \le \int_{\Rd}  h(t,x) dx \le 1.
	\end{equation*}
\end{proof}

\begin{theorem} 
	\label{th:convcons_MF}
	Suppose that Assumptions \ref{ass:well-posmicro_F} and \ref{ass:confrand_F} hold. Suppose that the hypotheses of Proposition \ref{prop:convsol_MF} are fulfilled with $C$ the compact subset given by Assumption \ref{ass:confrand_F}. 
	Then, denoting by $ x^{\alpha, f}_t[h]$ and $x^{\alpha, \hat{f}_M(\bY(\cdot))}_t[h(\bY(\cdot))]$ the consensus points \eqref{eq2: xalphat_f_bar} and \eqref{eq2: xalphat_M_bar}, it holds that
	\begin{equation*}
		\left| \;	x^{\alpha, f}_t[h] - x^{\alpha, \hat{f}_M(\bY(\cdot))}_t[h(\bY(\cdot))] \; \right| \to 0 \quad \tn{for $\alpha,t,M$ sufficiently large and w.p.1.}
	\end{equation*}
\end{theorem}

\begin{proof}[Proof of Theorem \ref{th:convcons_MF}]
	Let $C$ be the compact subset given by Assumption \ref{ass:confrand_F}.
	As in the proof of Proposition \ref{prop:convsol_MF}, let $\N \in \A$ be the set of measure zero (with respect to probability measure $\PP$) that guarantees the uniform convergence of $\hat{f}_M$ to $f$; let us fix $\omega^u \in \Omega \setminus \N$ and $\by^u := \bY(\omega^u)$. By means of a triangular inequality, we have that
	\begin{align*}
		\bigg| \; x^{\alpha, f}_t[h] &- x^{\alpha, \hat{f}_M(\by^u)}_t[h(\by^u)] \; \bigg| \\
		&\le 
		\underbrace{\left| \; x^{\alpha, f}_t[h] - x^{\alpha, \hat{f}_M(\by^u)}_t[h] \; \right|}_{(A)} + 
		\underbrace{\left| x^{\alpha, \hat{f}_M(\by^u)}_t[h] - x^{\alpha, \hat{f}_M(\by^u)}_t[h(\by^u)  \; \right|}_{(B)}.
	\end{align*}
	We now prove that $(A)$ and $(B)$ approach zero for $\alpha,t,M$ sufficiently large, hence obtaining the thesis.\\
	The numerator of $(A)$ reads $ \left| \int_{\Rd} x \left( e^{-\alpha f(x)} - e^{-\alpha \hat{f}_M(x, \by^u)} \right)  h(t,x) dx \right|$, 
	while the denominator is $ \left| \int_{\Rd} \left( e^{-\alpha f(x)} - e^{-\alpha \hat{f}_M(x, \by^u)} \right)  h(t,x) dx \right|$.
	By applying estimates \eqref{eq:estimatenum_diffobj} and \eqref{eq:estimateden_diffobj} of Lemma \ref{lem:convcons_MF} to the numerator and denominator resp. and by using that 
	\[
	\sup_{x \in C} \left| \; f(x) - \hat{f}_M(x, \by^u) \; \right| \xrightarrow[]{M \to \infty} 0,
	\]
	for how $\omega^u$ and $\by^u$ are chosen, we conclude the convergence of $(A)$ to $0$. 	
	Then, $(B)$ is estimated with \eqref{eq:estimatexalpha_diffmeas} of Lemma \ref{lem:convcons_MF} for $\by = \by^u$. As $C_0$ of Lemma \ref{lem:convcons_MF} is independent of $M$ and $\by^u$ (and clearly of $\alpha,t$), it possible to pass to the limit for $ M, \alpha, t$ large and use Proposition \ref{prop:convsol_MF} to conclude that $(B)$ converges to zero.
\end{proof}

Finally, we can complete the diagram of Figure \ref{diag:onlyN} in Figure \ref{diag:completeMF} with a vertical arrow between the mean-field formulations. The diagram is introduced only for visualization purposes and the vertical arrow expresses the two convergence results rigorously proven in Proposition \ref{prop:convsol_MF} and Theorem \ref{th:convcons_MF}. 
\begin{figure}[h!tbp]
	\centering
	\resizebox{13cm}{2.7cm}{%
		\begin{tikzpicture}[node distance=2cm]	
			\node (cboM) [CBOalg] {CBO algorithm \eqref{eq2: complete cbo for fhatM}\\ with consensus $x^{\alpha,\hat{f}_M(\bY(\cdot))}_t$ \eqref{eq2: xalphat_M}};
			\node (cbof) [CBOalg, above of=cboM] {CBO algorithm \eqref{eq2: complete cbo for f, limMinfty done}\\ with consensus $x^{\alpha, f}_t$ \eqref{eq2: xalphat_f}};
			\node (MFM) [MF, right of=cboM, xshift=6.5cm] {Mean-field formulation \eqref{eq2: mf eq complete for fM}\\with consensus $x^{\alpha, \hat{f}_M(\bY(\cdot))}_t[h(\bY(\cdot))]$ \eqref{eq2: xalphat_M_bar}};
			\node (MFf) [MF, above of=MFM] {Mean-field formulation \eqref{eq2: mf eq complete for EF}\\with consensus $x^{\alpha, f}_t[h]$ \eqref{eq2: xalphat_f_bar}};
			\draw [arrow,dashed,very thick] (cboM) -- node[anchor=north]{$N \to \infty$} (MFM);
			\draw [arrow,very thick,<->] (MFM) -- node[anchor=west]{$M, \alpha, t$ sufficiently big} node[anchor=east]{w.p.1} (MFf);
			\draw [arrow,dashed,very thick] (cbof) -- node[anchor=north]{$N \to \infty$}(MFf);
	\end{tikzpicture}}
	\caption{Derivation of two relations between mean-field formulations \eqref{eq2: mf eq complete for EF} and \eqref{eq2: mf eq complete for fM} and  w.p.1 and for $M,\alpha,t$ sufficiently large. The precise sense in which the convergence of the vertical arrow holds is given in Proposition \ref{prop:convsol_MF} and Theorem \ref{th:convcons_MF}. The sense in which the horizontal arrows holds is explained in Figure \ref{diag:onlyN}.}
	\label{diag:completeMF}
\end{figure}

\begin{remark}
	In Subsection \ref{subsec:SAA_MFlevelder}, we verified the well-posedness of algorithms  \eqref{eq2: complete cbo for f, limMinfty done} and \eqref{eq2: complete cbo for fhatM} under Assumption \ref{ass:well-posmicro_F}. We remark that requiring $x \mapsto F(x,\bfy)$ locally Lipschitz continuous for any $\bfy \in E$ would have been sufficient to guarantee the existence and uniqueness of solutions to the algorithms.	
	Analogously, if condition 2. of Assumption \ref{ass:well-posMF_F} held for constants $J_F(\bfy), c_u(\bfy) >0$ dependent on $\bfy \in E$, it would have been enough to imply the well-posedness of mean-field formulations \eqref{eq2: mf eq complete for fM} and \eqref{eq2: mf eq complete for EF}.
	However, the stronger conditions guaranteed the independence of $C_0(\alpha,J_F,K_4(L_F))$ and $C_1(K_4(L_F))$ of Lemma \ref{lem:convcons_MF} of $\by^u$, fact that was of paramount importance for taking the limit for $M \to \infty$ in Theorem \ref{th:convcons_MF}.
\end{remark}
\begin{remark}
	Throughout the section we assumed that the diffusion term of the CBO algorithm was of isotropic type. All the results may be extended to the anisotropic diffusion, provided that \rev{\cite[Theorem 1]{fornasier2022convergence}, \cite[Theorem 3]{fornasier2022convergence}} are used in place of \rev{\cite[Theorems 3.1,3.2]{carrillo2018analytical}, \cite[Theorem 3.7]{fornasier2021consensus}} respectively (and the assumption $2\lambda>d\sigma^2$ is relaxed to $2\lambda>\sigma^2$).
\end{remark}

\begin{remark}
	A question that arises \rev{spontaneously} at this point is whether the newly obtained diagram of Figure \ref{diag:completeMF} can be completed also with a vertical arrow at the microscopic level, namely if we can define the analogous of Proposition \ref{prop:convsol_MF} and Theorem \ref{th:convcons_MF} for CBO algorithms \eqref{eq2: complete cbo for f, limMinfty done} and \eqref{eq2: complete cbo for fhatM}.
	We are of the opinion that the answer is positive and in this remark we give an intuition as to why this is so.
	
	Let us begin with the discussion of the extension of Proposition \ref{prop:convsol_MF}.
	We take a step back and look at the two approaches available in the literature to prove the convergence of CBO-type algorithms. 
	The first technique, hereafter called ``variance-based'', was proposed in \cite{pinnau2017consensus,carrillo2018analytical} and further applied in many other works to prove convergence at the mean-field level \cite{pinnau2017consensus,carrillo2018analytical,carrillo2021consensus} and at the microscopic level \cite{ha2020convergence,ha2021convergence,ko2022convergence}. It consists of two steps: the investigation of the emergence of consensus from evolving the system variance and the assessment that the objective evaluated in the consensus is arbitrarily close to the optimal point. 
	The second, introduced in \cite{fornasier2021consensus,fornasier2022convergence} and hereafter called ``gradient flow-based'', is valid at the mean-field level and is based on the intuition that the mean-field is equivalent to a gradient flow with respect to the mean squared error. At variance with the previous technique, it proves the convergence in one step by showing the decay of a suitable functional. 
	Theorem 3.7 of \cite{fornasier2021consensus} used in Proposition \ref{prop:convsol_MF} exploits the just-mentioned strategy, but, to the best of our knowledge, an analogous of such theorem is missing at the algorithmic level. As a consequence, completing the diagram of Figure \ref{diag:completeMF} would require either to replicate the ``gradient flow-based'' approach at the microscopic level or to use the ``variance-based'' approach through \cite{ha2020convergence,ha2021convergence,ko2022convergence}. We remark that the presence of the double step in the second option prevents the use of direct estimation such as that of \eqref{eq:proofprop_W2triang}. 
	
	Subsequently, Theorem \ref{th:convcons_MF}'s proof is constructed on Proposition \ref{prop:convsol_MF} and on Lemma \ref{lem:convcons_MF}. The latter is immediately verified to hold true for the empirical distributions
	\begin{equation*}
		h^N(t,x) = \frac{1}{N} \sum_{i=1}^N \delta(X^i_t-x), \quad h^N(t,x)(\bY(\cdot)) = \frac{1}{N} \sum_{i=1}^N \delta(X^i_t(\bY(\cdot))-x)
	\end{equation*}
	associated to the solutions $\{ X^i_t\}^i_t$ of \eqref{eq2: complete cbo for f, limMinfty done} and $\{ X^i_t(\bY(\cdot))\}^i_t$ of \eqref{eq2: complete cbo for fhatM} resp. (and $\delta$ the Dirac delta function), and thus it is easily adapted to the microscopic scenario.
	We only observe that proving it would require some additional technical details with respect to Lemma \ref{lem:convcons_MF} as $h^N(t,x)(\bY(\cdot))dx$ is a random measure.
\end{remark}

\section{The quadrature approach}
\label{sec:quadr}

In Section \ref{sec:SAA}, we presented the sample average approximation approach and estimated $f(x) = \EE_\PP[F(x,\mathbf{Y})]$ with $\hat{f}_M(x,\bY)$ for any $x \in \Rd$ and $\bY$ $M$-dimensional sample from the random vector $\mathbf{Y}$.
We had two key parameters: $M$, the sample size, and $N$, the number of agents interacting in the search space $\Rd$; we derived the diagram in Figure \ref{diag:completeMF} by analyzing the limits in the two parameters independently.
In this section, we are interested in studying a ``joint limit'', namely to consider a discretization of $f$ that depends on a parameter which is set equal to the number of agents of the algorithm. Apart from the theoretical analysis of the equation obtained by letting the joint limit approach infinity (see Subsection \ref{subsec:MFquad}), which we believe is interesting in its own right, 
considering such limit allows us to quantify the gap between the microscopic and the mean-field level with one unique metric (see Subsections \ref{subsecnum: rate quadr} and \ref{subsecnum: higher k}).
In the following we assume that $\tn{supp}(\mathbf{Y}) = E = [e_l,e_u]^k \subset \R^k$ and require that $\nu^\mathbf{Y}$ is absolutely continuous with respect to the Lebesgue measure and call $\theta_\mathbf{Y}$ its density. For \rev{convenience of notation}, we define $\mathcal{E} = [e_l,e_u]$ and $l(\mathcal{E}) = e_u-e_l$.
We approximate 
\begin{equation*}
	f(x) = \int_{E} F(x,\mathbf{y}) \theta_\mathbf{Y}(\mathbf{y}) d\mathbf{y} = \int_{\mathcal{E}} \ldots \int_{\mathcal{E}} F(x,(y_1, \ldots, y_k)^T) \theta_\mathbf{Y}((y_1, \ldots, y_k)^T) dy_1 \ldots dy_k
\end{equation*}
by a composite midpoint quadrature formula (see e.g. \cite{hamming2012numerical}) with $Q \in \mathbb{N}$ points for each coordinate direction. 
Given the set of nodes
\begin{equation}
	\label{def3: discr points midpoint, multik}
	y_l^r = e_l + \frac{l(\mathcal{E})}{2Q}(2q-1), \quad r = 1, \ldots, Q, \; l = 1, \ldots, k,
\end{equation}
the formula reads
\begin{equation*}
	\tilde{f}_{Q^k}(x,\{y^r_l\}^{r = 1, \ldots, Q}_{l=1,\ldots,k}) = \left( \frac{l(\mathcal{E})}{Q} \right)^k \sum_{r_1=1}^{Q} \ldots \sum_{r_k=1}^{Q} F(x,(y^{r_1}_1,\ldots,y^{r_k}_k)^T)\theta_\mathbf{Y}((y^{r_1}_1,\ldots,y^{r_k}_k)^T).
\end{equation*}
We relabel the nodes so to have one running index $j = 1,\ldots,Q^k$ and obtain  
\begin{equation}
	\label{def3:set nodes}
	\{ \mathbf{y}^j = ((y_1,\ldots,y_k)^T)^j \}^{j=1,\ldots,Q^k}
\end{equation}
(a possible way to do the labeling consists of ordering in a matrix structure all the $Q^k$ nodes and then indexing them starting from the lower left corner and going up the rows). For consistency with the SAA approach, we introduce the notation $\by := (\mathbf{y}^1, \ldots, \mathbf{y}^{Q^k})$: then, the quadrature formula reduces to 
\begin{equation}
	\label{def3: ftildeQk}
	\tilde{f}_{Q^k}(x,\by) := \left( \frac{l(\mathcal{E})}{Q}\right)^k \sum_{j=1}^{Q^k} F(x,\mathbf{y}^j) \theta_\mathbf{Y}(\mathbf{y}^j).
\end{equation}
The approximation $f \approx \tilde{f}_{Q^k}$ leads to the optimization problem 
\begin{equation}
	\label{eq3: min problem det, quad}
	\min_{x \in \Rd} \tilde{f}_{Q^k}(x,\by),
\end{equation}
\rev{which we assume, similarly to \eqref{eqi: min problem main} and \eqref{eq2: min problem, saa}, has a unique minimizer, and}
to which we apply a CBO-type algorithm. 
Now, we choose the number of particles in the algorithm to be equal to the number of discretization points $Q^k$. 
Consequently, each agent $i$ has phase space values
\begin{equation}
	\label{def3: (Xit,yi)}
	(X^i_t,\mathbf{y}^i) \in \Rd \times E. 
\end{equation}
The dynamics of the agents takes place in the augmented space $\Rd \times E$.
Then, in addition to the system of SDEs for the position vectors, we have the system of SDEs for each node $\bfY^i_t$
\begin{equation*}
	d\mathbf{Y}^i_t =0 \quad \textnormal{for} \; i=1,\ldots,N \quad \tn{with initial condition $(\bfY^1_0,\ldots,\bfY^N_0) = \by$.}
\end{equation*}
Clearly, $(\bfY^1_t,\ldots,\bfY^N_t) = \by$ for any $t$. Finally, the CBO algorithm is given by 
\begin{subequations}
	\label{eq3: complete cbo for ftildeN}
	\begin{align}
		dX^i_t(\by) &= -\lambda(X^i_t(\by) - x^{\alpha,\tilde{f}_N(\by)}_t) dt + \sigma D^i_t(\by) dB^i_t \quad \textnormal{for} \; i=1,\ldots,N, \label{eq3: sde cbo for ftildeN}\\
		x^{\alpha,\tilde{f}_N(\by)}_t &= \frac{\sum_{i=1}^{N} X^i_t(\by) \exp(-\alpha \tilde{f}_N(X^i_t(\by),\by))}{\sum_{i=1}^{N} \exp(-\alpha \tilde{f}_N(X^i_t(\by),\by))}, \label{def3: xalphat_N}\\
		d\mathbf{Y}^i_t &=0 \quad \textnormal{for} \; i=1,\ldots,N \quad \iff \quad (\bfY^1_t,\ldots,\bfY^N_t) = \by \label{def:by_quad}
	\end{align}
\end{subequations}
where $\lambda, \sigma >0$, $\alpha >0$, $B^i_t$ are $d$-dimensional independent Brownian processes and $D^i_t(\by)$ is equal to
\begin{equation*}
	D^i_{t,\tn{iso}}(\by) := |X^i_t(\by) - x^{\alpha,\tilde{f}_N(\by)}_t| I_d 
	\quad \tn{or} \quad  D^i_{t,\tn{aniso}}(\by) := \tn{diag}(X^i_t(\by) - x^{\alpha,\tilde{f}_N(\by)}_t).
\end{equation*}
We require the system to be supplemented with the same initial conditions $\{X^i_0\}^i$ of Section \ref{sec:SAA} and, for the sake of simplicity, that $D^i_t(\by) = D^i_{t,\tn{iso}}(\by)$.\\

Microscopic system \eqref{eq3: complete cbo for ftildeN} is well-posed, provided that $ \EE [|(X^1_0,\ldots,X^N_0)|^2] < \infty $.
Indeed, the assumptions on $(X^1_0,\ldots,X^N_0)$ imply that $((X^1_0,\bfY^1_0)^T,\ldots,(X^N_0,\bfY^N_0)^T)$ has finite-second moment. Then, $\tilde{f}_N$ is lower bounded and it is straightforward to verify that Assumption \ref{ass:well-posmicro_F} implies that $x \mapsto \tilde{f}_N(x,\by)$ is locally Lipschitz continuous with constant $L_F \frac{1}{N} \sum_{i=1}^N \theta_{\bfY}(\bfy^j)$. All in all, the existence of a unique strong solution 
$\{X^i_t(\by)\}^i_t$ is guaranteed by \cite[Theorem 2.1]{carrillo2018analytical}, \cite[Theorem 3.1]{durrett1996stochastic}.

\subsection{The mean-field equation for $N \to \infty$}
\label{subsec:MFquad}

In this section we obtain the mean-field formulation associated to system \eqref{eq3: complete cbo for ftildeN} in the regime of the large number of particles.
More precisely, let $G^N(t)$ be the $N$-particle probability distribution over $(\Rd \times E)^N$ at time $t$. Then, the propagation of chaos assumption on the marginals 
\footnote{\rev{
		At the beginning of Section 3, we introduce $\tilde{f}_{Q^k}$ and explain that it can be rewritten to feature a single running index. 
		We note that there is no unique way to accomplish this rewriting: this is the key point in understanding how the problem setup aligns with typical propagation of chaos framework. We fix the nodes at the positions determined by the quadrature formula, then sample the indices and assign them to the points that we use to initialize the algorithm ($\by$ appearing in \eqref{def:by_quad}). In this way, the initial data satisfies a chaoticity condition.}}
translates to the approximation $G^N(t) \approx g(t)^{\otimes N}$ with $g(t)$ being probability distribution over $\Rd \times E$.
Let 
\begin{equation}
	\label{def3: w0}
	w_0(\mathbf{y}) = \left( \frac{1}{l(\mathcal{E})} \right)^k, \quad \tn{for $\mathbf{y} \in E$},
\end{equation}
be the probability density reconstructed with the nodes $\left\lbrace  \mathbf{y}^j \right\rbrace^{j=1,\ldots,N}$ specified by \eqref{def3:set nodes}. Then, it is straightforward to observe that $g$ has moments
\begin{equation*}
	\int_{\Rd} g(t,x,\mathbf{y}) dx = w_0(\mathbf{y}), \quad \rho(t,x) := \int_{E} g(t,x,\mathbf{y}) d\mathbf{y}
\end{equation*}
and that it \rev{formally} solves in the distributional sense the non-linear Fokker-Planck equation 
\begin{subequations}
	\label{eq3: mf eq complete quad}
	\begin{align}
		\pd_t g(t,x,\mathbf{y}) &= \lambda \nabla \cdot \left(  (x - x^{\alpha,f}_t[\rho] ) \:g(t,x,\mathbf{y}) \right) 
		+\frac{\sigma^2}{2} 
		\Delta_x \left( |x-x^{\alpha,f}_t[\rho]|^2 g(t,x,\mathbf{y}) \right)
		\label{eq3: mf eq quad}\\
		x^{\alpha,f}_t[\rho] &:= \frac{\int_{\Rd} x \exp(-\alpha f(x)) \rho(t,x) dx }{\int_{\Rd} \exp(-\alpha f(x)) \rho(t,x) dx}, \label{def3: xalphat_N, lim}\\
		w_0(\mathbf{y}) &\mu_0(dx) d\mathbf{y} = \lim_{t \to 0} g(t,x,\mathbf{y})  dx d\mathbf{y}. \label{eq3: in cond mf eq quad}
	\end{align}
\end{subequations}
The independence of the consensus point \eqref{def3: xalphat_N, lim} of $\bfy$ is expected as a consequence of the fact that the moment of $g$ with respect to $x$ is equal to $w_0$.

\begin{remark}
	\label{rem:equalityMF}
	If $h(t,x)dx$ is a weak solution to \eqref{eq2: mf eq complete for EF}, namely the mean-field formulation for the true objective function $f$, then $h(t,x)dx$ satisfies \eqref{eq3: mf eq complete quad} in the distributional sense.
	Indeed, it is sufficient to choose $g(t,x,\bfy) = h(t,x)$ 
	in the weak formulation of \eqref{eq3: mf eq complete quad} and observe that the equation obtained is exactly the weak formulation of \eqref{eq2: mf eq complete for EF}. Using then that $h$ is a solution of \eqref{eq2: mf eq complete for EF}, we obtain the desired result. 
	In addition, thanks to the uniqueness of the solution to \eqref{eq3: mf eq complete quad} (which is guaranteed if we derived the mean-field formulation rigorously using \cite{huang2022mean}), also the stronger statement that all weak solutions to \eqref{eq2: mf eq complete for EF} satisfy \eqref{eq3: mf eq complete quad} in the distributional sense holds. 
	The just mentioned result may be visualized as follows:
	\begin{figure}[h!tbp]
		\centering
		\resizebox{13cm}{3.5cm}{%
			\begin{tikzpicture}[node distance=1.5cm]
				\node (cboN) at (0,0) [CBOalg] {CBO algorithm \eqref{eq3: complete cbo for ftildeN}\\ with consensus point $x^{\alpha,\tilde{f}_N(\by)}_t$ \eqref{def3: xalphat_N}};
				\node (MFf2) at (7.5,1.3) [MF2] {Mean-field formulation \eqref{eq3: mf eq complete quad}\\with consensus point $x^{\alpha,f}_t[\rho]$ \eqref{def3: xalphat_N, lim}};
				\node (MFf) [MF, above of=MFf2] {Mean-field formulation \eqref{eq2: mf eq complete for EF}\\with consensus point $x^{\alpha, f}_t[h]$ \eqref{eq2: xalphat_f_bar}};
				\draw  (MFf2) -- (MFf);
				\draw [arrow,blue,very thick,dashdotted] (cboN) -- node[anchor=east, yshift=-0.5cm,xshift=1cm]{$N \to \infty$}(MFf2);
		\end{tikzpicture}}
		\caption{Derivation of the relations between CBO algorithm \eqref{eq3: complete cbo for ftildeN} and mean-field formulations \eqref{eq2: mf eq complete for EF} and \eqref{eq3: mf eq complete quad}. The limit in $N$ holds under the propagation of chaos assumption on the marginals; the vertical arrow connecting the two mean-fields should be intended in the sense expressed by Remark \ref{rem:equalityMF}.}
		\label{diagi: diagSAA_completed}
	\end{figure}
\end{remark}

\section{Numerical experiments}
\label{sec:numerics}

In this section we validate the two proposed approaches by some numerical tests. 
We investigate their rates in suitable norms and for various values of $k$, dimension of the random space (Tests $1,2,3$ of Subsections \ref{subsecnum: rate SAA}, \ref{subsecnum: rate quadr}, \ref{subsecnum: higher k}). This allows us to equip the diagrams of Figures \ref{diag:completeMF} and \ref{diagi: diagSAA_completed} with the rates as $N,M \to \infty$ and $N \to \infty$ respectively.
Thereafter, we compare the success rates of the two methods on an example drawn from the field of economics (Test $4$ of Subsection \ref{subsecnum:test4_succrate}). \\

The idea common to both procedures is the approximation of $f$ by a discretization, either $\hat{f}_M$ \eqref{def2: fhatM} or $\tilde{f}_N$ \eqref{def3: ftildeQk}, and then the resolution of the newly obtained problem by a CBO algorithm.
We adopt an explicit Euler-Maruyama scheme, see e.g. \cite{higham2001algorithmic}, to numerically simulate the solutions to the systems of SDEs \eqref{eq2: sde cbo for fhatM} and \eqref{eq3: sde cbo for ftildeN}.
We discretize the time interval $[0,T]$ by introducing the time step $\Delta t = \frac{T}{n_{\tn{it}}-1}$ and by defining the $n_{\tn{it}}$ points $t_{h} = ({h}-1)\Delta t, {h} = 1, \ldots, n_{\tn{it}}$.
Given the $i$th particle, we set $X^i_{h}$ as the approximation of $X^i_{t_{h}}$, so that the scheme reads
\begin{equation}
	\label{eq5: discretized cbo for M,N}
	X^i_{{h}+1} = X^i_{{h}}-\lambda(X^i_{h} - x^{\alpha, \#}_{h}) \Delta t + \sigma D^i_{h} \sqrt{\Delta t} Z^i_h \quad \textnormal{for} \; {h}=1,\ldots,n_{\tn{it}}-1,
\end{equation}
where the ${h}$th increment $dB^i_{h}$ of the Brownian motion associated to agent $i$ is approximated by $\sqrt{\Delta t} Z^i_h$, with $Z^i_h$ normally distributed with zero mean and identity covariance, and where $x^{\alpha, \#}_{h}$ stands for either $x^{\alpha,\hat{f}_M(\bY(\cdot))}_{t_h}$ \eqref{eq2: xalphat_M} or $x^{\alpha,\tilde{f}_N(\by)}_{t_h}$ \eqref{def3: xalphat_N}. 
There are two sources of stochasticity in the scheme defined above: the probability distribution $\mu_0$, from which the initial positions $X^1_0, \ldots, X^N_0$ are sampled, and the ${h}$th increment $dB^i_{h}$ of the Brownian motion. In our numerical tests we run the algorithm $n_{samplesCBO}$ times to accommodate for such uncertainties. 

We simulate the solutions to the mean-field formulations presented in the work by choosing a large number of particles, hereafter denoted by $N_{ref}$. An alternative approach could consist, for instance, of simulating the solution to the Fokker-Planck equations with a splitting scheme in time and a discontinuous Galerkin method in space as done in \cite{pinnau2017consensus}.

\begin{remark}
	Let us fix an iteration $h$ in the update rule \eqref{eq5: discretized cbo for M,N}. The computational complexity of evaluating $X^i_{h+1}$ is governed by the number of evaluations of $F$ and by the $O(N)$ computations required to calculate $x^{\alpha, \#}_{h}$. One can reduce the computational complexity of the calculation of $x^{\alpha, \#}_{h}$ by considering mini-random batch techniques \cite{albi2013binary,jin2020random}:
	they select a random subset $I^R_h \subset \{1,\ldots,N\}$ of $R \ll N$ particles and substitute
	$\frac{1}{N} \sum_{i=1}^N (\cdot)^i$ with $\frac{1}{R} \sum_{i \in I^R_h}^R (\cdot)^i,$
	in the computation of $x^{\alpha,\#}_h$, so that the evaluation of the new consensus point lowers from $O(N)$ to $O(R)$.
\end{remark}
In the following, we restore the time variable $t \ge 0$ to indicate one of the nodes of the time mesh $\{ t_1, \ldots, t_{n_{\tn{it}}}\}$. 

\subsection{Test 1: Convergence rates for SAA and CBO}
\label{subsecnum: rate SAA}

\begin{figure}[h!tbp]
	\centering
	\subfigure[]{\includegraphics[scale = 0.28]{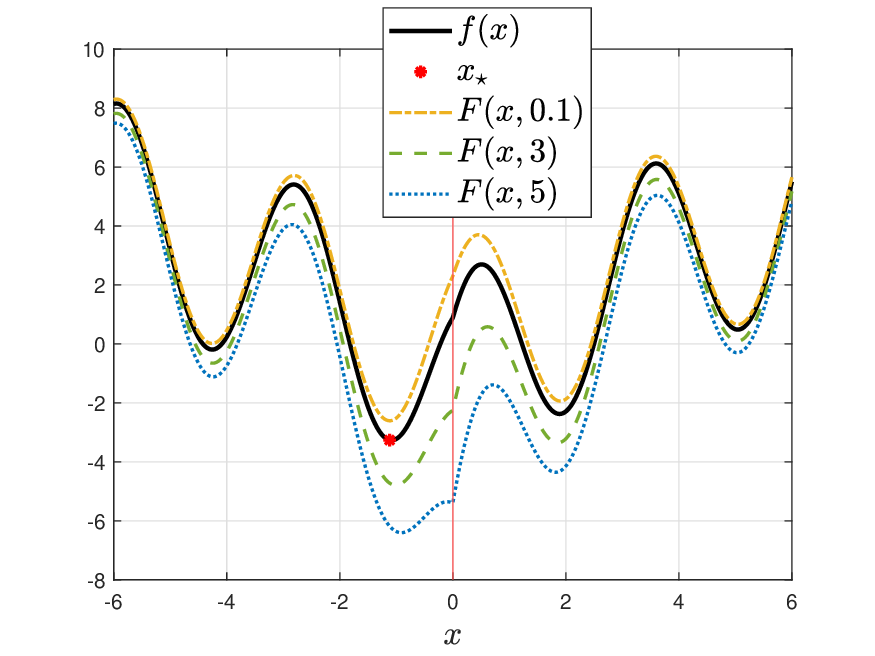}} 
	\subfigure[]{\includegraphics[scale = 0.28]{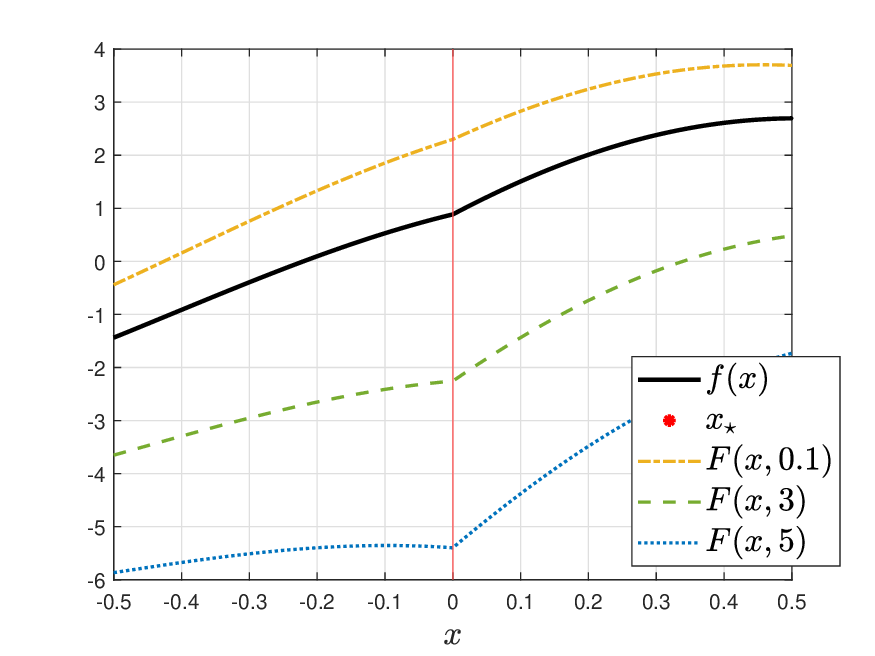}} 
	\subfigure[]{\includegraphics[scale = 0.28]{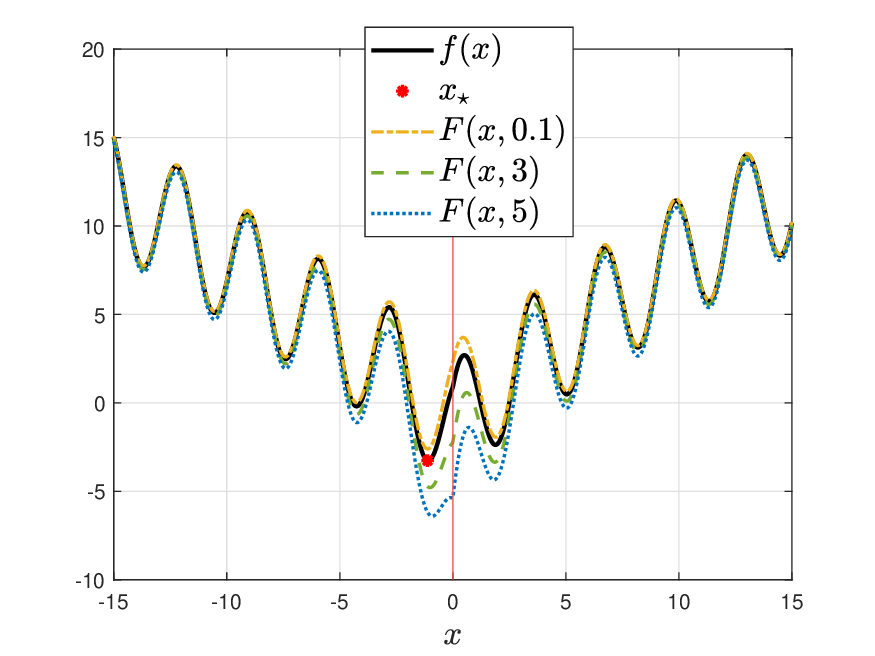}}
	\caption{Plot of $F$ \eqref{def: ackleystocmodR1F} and $f$ \eqref{def: ackleystocmodR1f} for the choices $Y = 0.1,3,5$. In plot (a), $x \in [-6,6]$, in plot(b), $x \in [-0.5,0.5]$, in plot(c), $x \in [-15,15]$. Plot (b) aims at highlighting the non-differentiability of $f(\cdot)$ and $F(\cdot,Y)$, for any $Y \in \R$, at $x=0$: we remark that the lost of regularity doesn't influence our CBO algorithms and theoretical results, as they all hold in a non-differentiability setting. Plot (c) shows that, for $x$ sufficiently large, $x \mapsto F(x,Y)$ is independent of $Y$.}
	\label{fig:ackleystocmodR1_plot&zoomx0}
\end{figure}

In this section we investigate the rates of the diagram of Figure \ref{diag:completeMF} of Section \ref{sec:SAA}. We design a problem in which the cost function admits a closed-form expression for the expected value and that satisfies all the assumptions of the section:

\begin{equation}
	\label{def: ackleystocmodR1F}
	F(x,Y) = e^{-0.2} (|x|+3(\cos(2x)+\sin(2x))) + Y \left(\arctan(|x|)-\frac{\pi}{2}\right), 
\end{equation}
for $x,Y \in \R$, and
\begin{equation}
	\label{def: ackleystocmodR1f}
	f(x) = \EE [F(x,Y)] = e^{-0.2} (|x|+3(\cos(2x)+\sin(2x))) + \EE[Y] \left(\arctan(|x|)-\frac{\pi}{2}\right), 
\end{equation}
for $x \in \R$.
A plot of $f$ and $F$ for three choices of $Y$ is given in Figure \ref{fig:ackleystocmodR1_plot&zoomx0}.
$f$ defined in \eqref{def: ackleystocmodR1f} is inspired by the cost function $4$ of the survey \cite{jamil2013literature} and it is continuous, multimodal/nonconvex and admits a unique minimizer at $\xs = -1.119$. The at-least-local Lipschitzianity of $x \mapsto |x|, \cos(2x), \sin(2x), \arctan(|x|)$ guarantees the satisfaction of $F$ \eqref{def: ackleystocmodR1F} of Assumptions \ref{ass:well-posmicro_F} and point 2. of Assumption \ref{ass:well-posMF_F}. In addition, $F$ is bounded from below and it grows quadratically in the far field, hence fulfilling the rest of Assumption \ref{ass:well-posMF_F}.
The term $Y \left(\arctan(|x|)-\frac{\pi}{2}\right)$ converges to 0 as $|x| \to \infty$, so that, for $x$ sufficiently large, $F$ is independent of $Y$, namely it satisfies Assumption \ref{ass:confrand_F} \footnote{We remark that the $\arctan$ term could be replaced by any function $u:\R \to \R$ such that $\lim_{|x|\to \infty} u(x) = 0$ (for instance, $u(x) = x e^{-0.1x^2}$).}.

\noindent In the following, we implement the isotropic CBO algorithm described by rule \eqref{eq5: discretized cbo for M,N} with the choices $d=1$ and
\begin{subequations}
	\label{param:tests1,2_complete}
	\begin{align}
		&Y \sim \mathcal{U}(\mathcal{E}), \quad \tn{with $\mathcal{E}= [0.1,1.9]$},  \label{param: ackleystocmodR1 Y} \\
		T = 10, \Delta t = 0.1;& \quad \lambda = 1, \alpha = 40, \sigma = 0.5; \quad \mu_0 = \mathcal{U}([-3,3]^d), \label{param5: test1,2}
	\end{align}
\end{subequations} 
where $\mathcal{U}$ denotes the continuous uniform distribution.
We remark that $\lambda, \sigma$ are chosen so to satisfy the inequality $2 \lambda > d\sigma^2$. 
Given the above parameters, it is straightforward to assess that $f$, $\mu_0$ and $\hat{f}_M$ (given by \eqref{def2: fhatM}) satisfy the remaining hypotheses of Proposition \ref{prop:convsol_MF} and Theorem \ref{th:convcons_MF}.\\

We fix $N = 5000$ and look at the gap between mean-field formulations \eqref{eq2: mf eq complete for EF} and \eqref{eq2: mf eq complete for fM}. Proposition \ref{prop:convsol_MF} and Theorem \ref{th:convcons_MF} give us a certain assurance (results w.p.1) that the error between the solutions and the consensus points approaches zero in the limit as the sample size $M$ grows to infinity (and for $\alpha,t$ sufficiently big). However, they do not provide any insight on the magnitude of the error for a given $M$. In this subsection, we investigate numerically such rate of convergence.\\
The estimation of the rates with which the optimal point $\hat{x}_M(\bY(\cdot))$ and the optimal value $\hat{f}_M(\hat{x}_M(\bY(\cdot)), \bY(\cdot))$ of the SAA problem \eqref{eq2: min problem, saa} converge to their correspective true counterpart $\xs$ and $f(\xs)$ constitutes a building block of the SAA theory together with the analysis of their consistency \cite{shapiro2003monte,shapiro2021lectures} (see also Remark \ref{rem:unifassumpSAA}). 
The consistency of the estimators is usually proven in the almost sure/w.p.1 sense (hence, our choice in Proposition \ref{prop:convsol_MF} and Theorem \ref{th:convcons_MF}), while the error usually considered in the rates of convergence is the root mean square or $L^2$ error (see also Section 2 of \cite{caflisch1998monte}). In the following, we proceed as in the standard theory and  carry out our numerical tests for 
\begin{equation}
	\label{def5: err_ratesaa, N big}
	\left( \EE_{\PP} \left[\left(  \bar{x}^{\alpha,\hat{f}_M(\bY(\cdot))}_t - \bar{x}^{\alpha,f}_t \right)^2 \right] \right)^{1/2},
\end{equation}
where, for brevity of notation, we denoted the mean-field regime with a bar sign rather than expliciting the mean-field measure ($x^{\alpha,f}_t[h] = \bar{x}^{\alpha,f}_t$ and other analogously).
To accommodate for the stochasticity induced by the Euler-Maruyama scheme, we run the CBO algorithm $n_{samplesCBO}=10$ times and we compute \eqref{def5: err_ratesaa, N big} numerically as
\begin{align*}
	&\left( \frac{1}{n_{{samplesY}}} \sum_{j=1}^{n_{{samplesY}}} \left( \frac{1}{n_{{samplesCBO}}} \sum_{u=1}^{n_{{samplesCBO}}} (x^{\alpha,\hat{f}_M(\bY(\omega^{(j)})), (u)}_t - x^{\alpha,f,(u)}_t) \right)^2 \right)^{1/2}
\end{align*}
where we have replaced $\EE_{\PP}$ with an average over $n_{samplesY}=200$.

\begin{figure}[h!tbp]
	\centering
	\subfigure[]{\includegraphics[scale = 0.4]{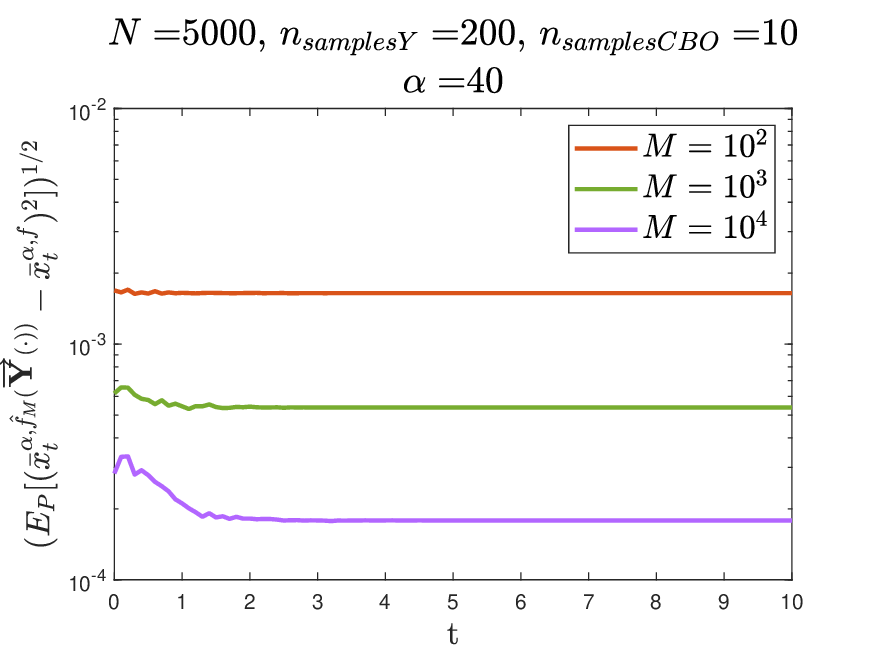}} 
	\subfigure[]{\includegraphics[scale = 0.4]{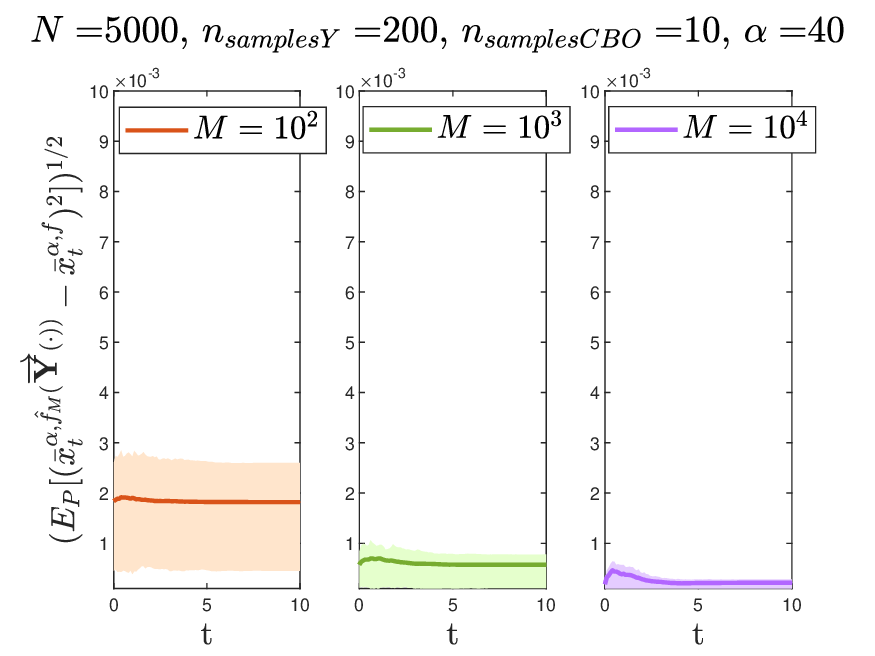}}  
	\subfigure[]{\includegraphics[scale = 0.4]{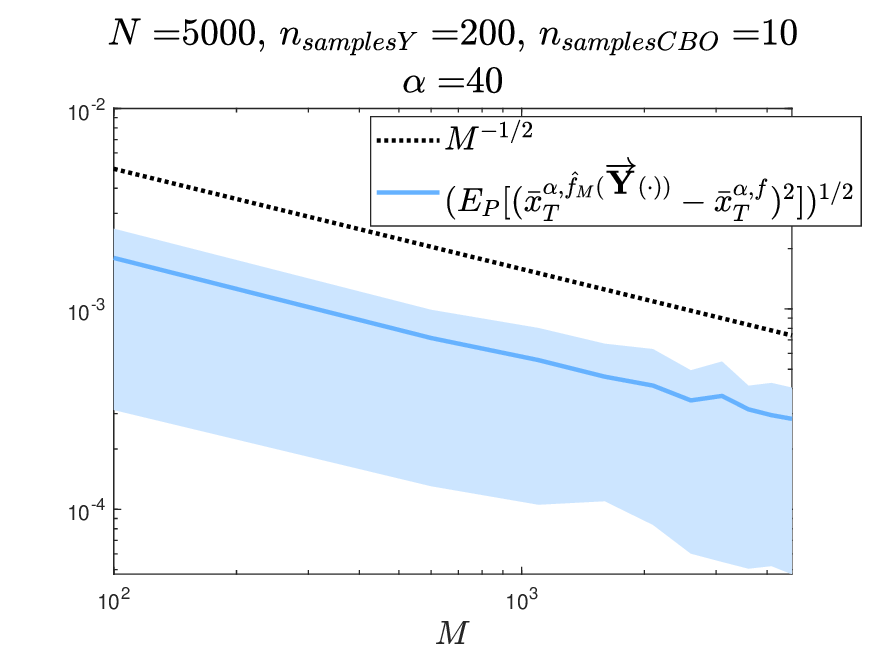}} 
	\caption{In the first row, evolution of the average error \eqref{def5: err_ratesaa, N big} (denoted by a continuous line) as a function of time for $M=10^2,10^3,10^4$. In the second row, evolution of the same error as a function of $M$ (varying from $10^2$ to $10^4$ with step $500$) and $t=T$. The shaded area of plots (b) and (c) indicates the $[0.15,0.85]$ quantile interval related to the uncertainty of repeatedly sampling of $\mathbf{Y}$. Plot (a) has a log-scale on $y$-axis, (b) a natural scale on $y$-axis and (c) a log scale on both axis.
		The objectives considered are \eqref{def: ackleystocmodR1F} and \eqref{def: ackleystocmodR1f} with \eqref{param: ackleystocmodR1 Y}, parameters are set to \eqref{param5: test1,2}.}
	\label{fig5:ackleystocmodR1_saacbo_saarateM}
\end{figure}

We present our numerical results in Figure \ref{fig5:ackleystocmodR1_saacbo_saarateM}.  

\noindent In the first row we proceed as in \cite{borghi2023kinetic} and investigate the evolution of error \eqref{def5: err_ratesaa, N big} as a function of time for three values of $M$. The second row  depicts the same error as a function of $M$ and for a fixed $t=T$. 
Plot (a) shows a decay in the average \eqref{def5: err_ratesaa, N big} for increasing $M$ and for each time $t$. 
This suggests a rate in $M$, which is quantitatively assessed in plot (c). We remark that the choice of the final time in such plot is fixed by the equidistance of the averages during the computation observed in plot (a). 
The rate $O(M^{-1/2})$ is coherent with the results of the classical theory of SAA, where such rate has been proven 
for the optimal points estimators provided that $F$ satisfies suitable differentiability and lipschitzianity assumptions around the minimizer \cite{shapiro2003monte,shapiro2021lectures}.
In addition, the fact that the averages (orange, green and purple continuous lines) are non-increasing nor non-decreasing suggests that the error introduced by approximating $\hat{f}_M$ with $f$ (and hence $x^{\alpha,\hat{f}_M(\bY(\cdot))}_t$ with $x^{\alpha,f}_t$) doesn't increase during the evolution of the particles' positions.
This observation is again coherent with the classical Monte Carlo theory \cite{caflisch1998monte,shapiro2021lectures,pareschi2013interacting}.
\begin{remark}
	\label{rem2: improvements saa}
	A possible disadvantage of Monte Carlo (MC)/SAA type techniques is the slow rate of convergence $O(M^{-1/2})$ (already mentioned in Remark \ref{rem:unifassumpSAA}. The exponent $1/2$ is intrinsically related to the choice of the root mean squared error and the MC approach and hence it can't be improved: however, the overall rate depends also on the variance of the function to be integrated. Several techniques have been proposed in literature to reduce such variance and hence to accelerate the converge rate (see e.g. \cite{caflisch1998monte,shapiro2021lectures}). 
\end{remark}
Next, we fix $M$ and let $N \to \infty$, as we are interested in the rate of convergence of the interacting particles systems \eqref{eq2: complete cbo for f, limMinfty done} and \eqref{eq2: complete cbo for fhatM} to the solutions to the mean-field equations \eqref{eq2: mf eq complete for EF} and \eqref{eq2: mf eq complete for fM} respectively (arrow connecting the pink and white boxes, first and second row resp.). 
This rate, known as mean-field approximation rate, has been theoretically investigated in the last decade: we mention \cite{fornasier2021consensus,fornasier2022convergence}
and \cite{gerber2023mean}, where the estimates are extended to stronger metrics than convergence in probability, in particular to Wasserstein distances.
We consider
\begin{equation}
	\label{def5: err_rateMF, M small}
	\EE_{\PP} \left[ W_p \left(\mu^{N,\hat{f}_M(\bY(\cdot))}_t, \bar{\mu}^{\hat{f}_M(\bY(\cdot))}_t \right) \right],
\end{equation}
for $p=1,2$, where we have added the expectation on $\PP$ with respect to \cite{gerber2023mean} to accommodate for the fact that the cost function $\hat{f}_M$ depends on the random $\mathbf{Y}$.
$\mu^{N,\hat{f}_M(\bY(\cdot))}_t, \bar{\mu}^{\hat{f}_M(\bY(\cdot))}_t$ denote the empirical measures associated to the samples \\
$\{X^i_t(\bY(\cdot))\}^{i=1,\ldots,N}$,
for a fixed $N$, and 
$\{X^i_t(\bY(\cdot))\}^{i=1,\ldots,N_{ref}}$,
for $N_{ref} = 10^5$, solutions to \eqref{eq2: complete cbo for fhatM}. 
We numerically simulate \eqref{def5: err_rateMF, M small} with
\begin{equation*}
	\hspace{-0.2cm}
	\frac{1}{n_{{samplesY}}} \sum_{j=1}^{n_{{samplesY}}} \left( \frac{1}{n_{{samplesCBO}}} \sum_{u=1}^{n_{{samplesCBO}}} W_p(\mu^{N, \hat{f}_M(\bY(\omega^{(j)})), (u)}_t, \bar{\mu}^{\hat{f}_M(\bY(\omega^{(l)})), (u)}_t ) \right)
\end{equation*}
and $n_{samplesY}=200$, $n_{samplesCBO} = 10$.

\begin{figure}[h!tbp]
	\centering
	\subfigure[]{\includegraphics[scale = 0.4]{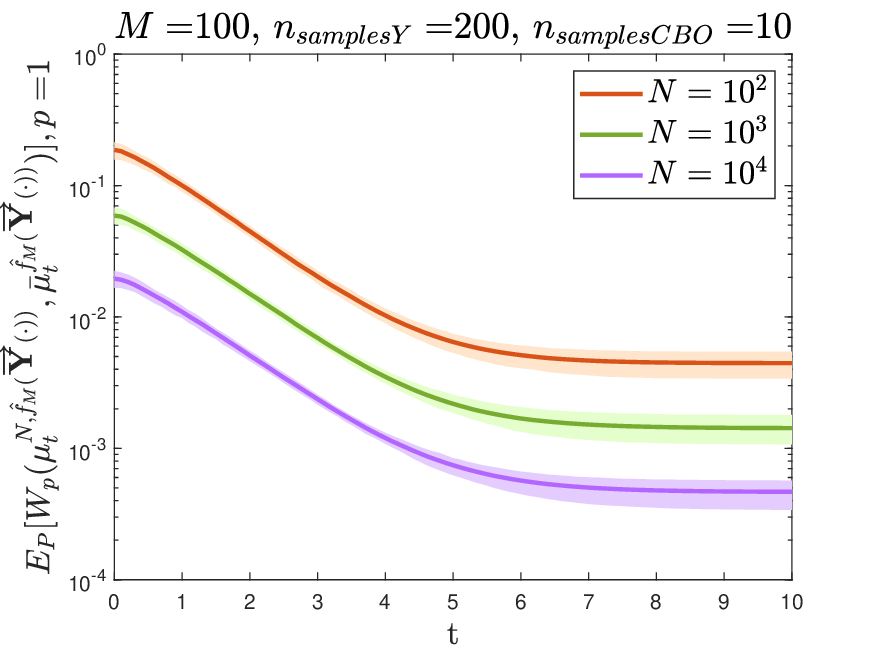}} 
	\subfigure[]{\includegraphics[scale = 0.4]{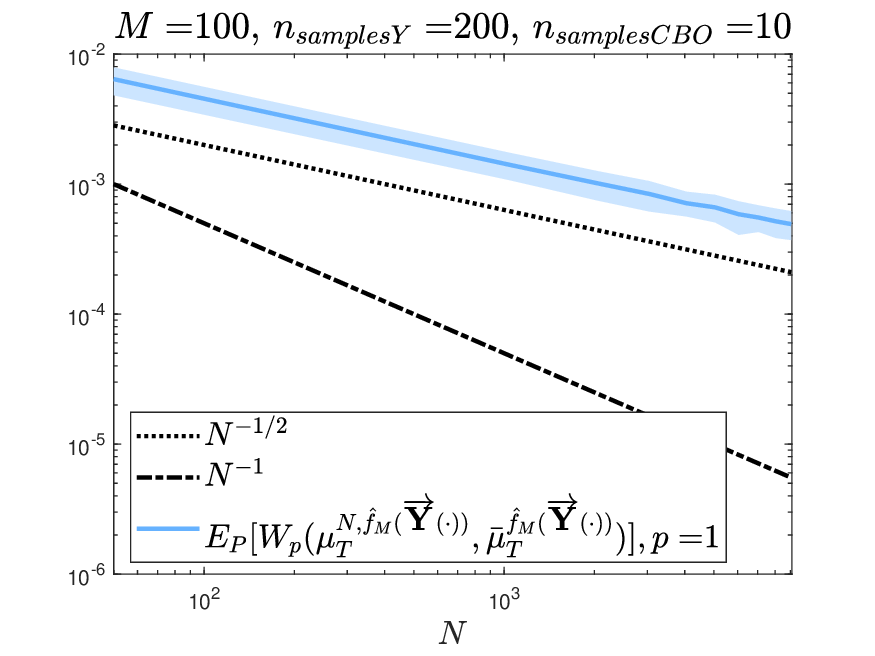}} 
	\subfigure[]{\includegraphics[scale = 0.4]{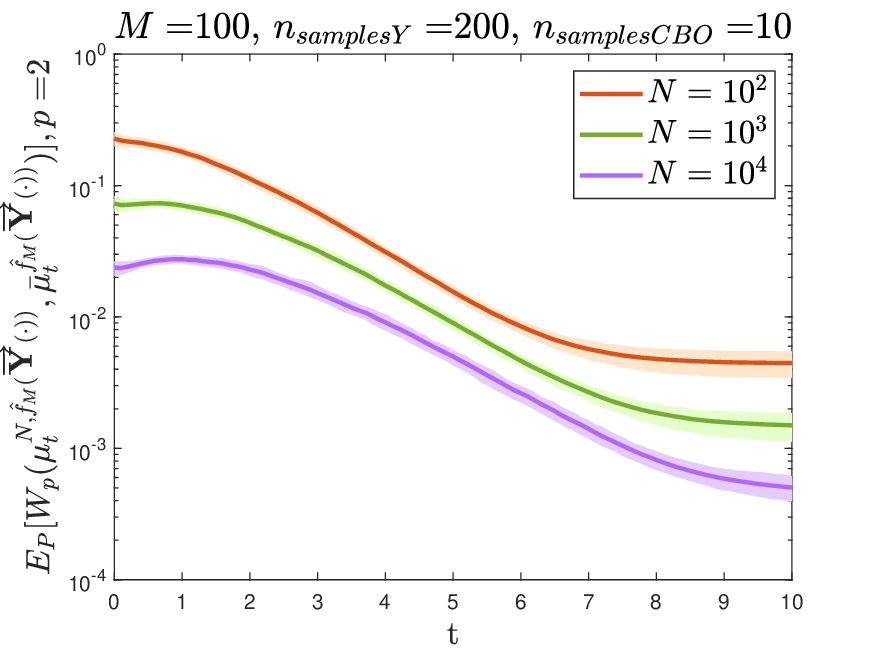}} 
	\subfigure[]{\includegraphics[scale = 0.4]{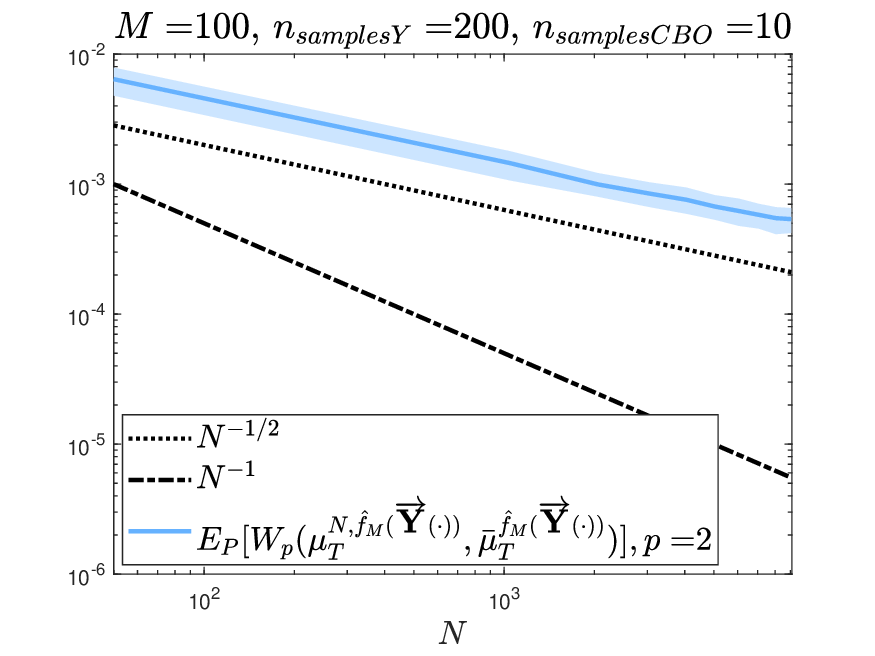}} 
	\caption{In the column on the left, evolution of the average error \eqref{def5: err_rateMF, M small} (denoted by a continuous line) as a function of time for $N=10^2,10^3,10^4$. In the column on the right, evolution of the same error as a function of $N$ (varying from $10^2$ to $10^4$ with step $500$) and $t=T$. 
		The first row corresponds to the $1$-Wasserstein distance, the second on the $2$-Wasserstein distance.
		The shaded area indicates the $[0.15,0.85]$ quantile interval related to the uncertainty of repeatedly sampling of $\mathbf{Y}$.
		Plot (a)-(c) have a log-scale on $y$-axis, while plot (b)-(d) a log scale on both axis.
		The objectives considered are \eqref{def: ackleystocmodR1F} and \eqref{def: ackleystocmodR1f} with \eqref{param: ackleystocmodR1 Y}, parameters are set to \eqref{param5: test1,2}.}
	\label{fig5:ackleystocmod_saacbo_MFrate2M}
\end{figure}

We present our main results in Figure \ref{fig5:ackleystocmod_saacbo_MFrate2M}, for $M=100$ and $p=1$ in the first row, $p=2$ in the second. 
The column on the left corresponds to error \eqref{def5: err_rateMF, M small} plotted as a function of time and for three values of $N$.
The right shows the same error plotted as a  function of $N$ and for fixed $t=T$. 
As in Figure \ref{fig5:ackleystocmodR1_saacbo_saarateM} (a)-(b), we observe, for a fixed time $t$ and for increasing $N$, a decrease in the average of the $p$th Wasserstein distance: numerically this validates the propagation of chaos assumption on the marginals; in addition, the decrease of the three averages as $t$ increases suggests a decrease in the mean-field error as the particles converge to the global minimizer $\xs$ \cite{borghi2023kinetic}. The rate exponent $1/2$ in Figure \ref{fig5:ackleystocmod_saacbo_MFrate2M} (b)-(d) is coherent with the analytical estimates presented in \cite{fornasier2021consensus,fornasier2022convergence,gerber2023mean}.
\begin{remark}
	Changing the value of $M$ in our simulations doesn't change the qualitative behavior of the plots in Figure \ref{fig5:ackleystocmod_saacbo_MFrate2M}. Indeed, $M$ governs the quality of the approximation $\hat{f}_M \approx f$, but, once it is fixed, $\hat{f}_M$ becomes an objective function to which to apply the classical results of the CBO on the mean-field approximation rate \cite{fornasier2021consensus,fornasier2022convergence,gerber2023mean}.
\end{remark}

\subsection{Test 2: Convergence rates for quadrature and CBO}
\label{subsecnum: rate quadr}

In Section \ref{sec:quadr}, we derived mean-field equation \eqref{eq3: mf eq complete quad} on the extended phase space associated to CBO algorithm \eqref{eq3: complete cbo for ftildeN} in the regime $N \to \infty$ and observed that it coincides with \eqref{eq2: mf eq complete for EF} (see Remark \ref{rem:equalityMF}). In this subsection we investigate the joint-limit error 
corresponding to the diagonal in the diagram of Figure \ref{diagi: diagSAA_completed}, that is
\begin{equation}
	\label{def:newrate_quad}
	W_p \left(\mu^{N,\tilde{f}_N(\by)}_t, \bar{\mu}^f_t \right) \approx \frac{1}{n_{{samplesCBO}}} \sum_{u=1}^{n_{{samplesCBO}}} W_p(\mu^{N, \tilde{f}_N(\by), (u)}_t, \bar{\mu}^{f, (u)}_t ),
\end{equation}
where ${\mu}^{N,\tilde{f}_{N}(\by)}_t$ denotes the empirical measure associated to $\{X^i_t(\by)\}^{i=1,\ldots,N}$, solution to
\eqref{eq3: complete cbo for ftildeN} for a fixed $N$
\footnote{\rev{
		When implementing microscopic CBO system (3.6), we fix the nodes $\by$ initially and apply the CBO update rule solely for the positions. This approach ensures that we do not increase the complexity of the algorithm (as we would if we also updated the nodes over time using the trivial dynamics $d\bfY^i_t = 0$) and highlights that the extension to augmented space of Section \ref{sec:quadr} is intended only for analytical purposes.}},
while $\bar{\mu}^f_t$ denotes the empirical measure associated to $\{X^i_t\}^{i=1,\ldots,N_{ref}}$, solution to \eqref{eq2: complete cbo for f, limMinfty done}  \footnote{We recall that in our numerical tests we simulate the mean-field equations with $N_{ref}$ particles and that \eqref{eq2: mf eq complete for EF} is obtained from \eqref{eq2: complete cbo for f, limMinfty done}.}.

\begin{figure}[h!tbp]
	\centering
	\subfigure[]{\includegraphics[scale = 0.45]{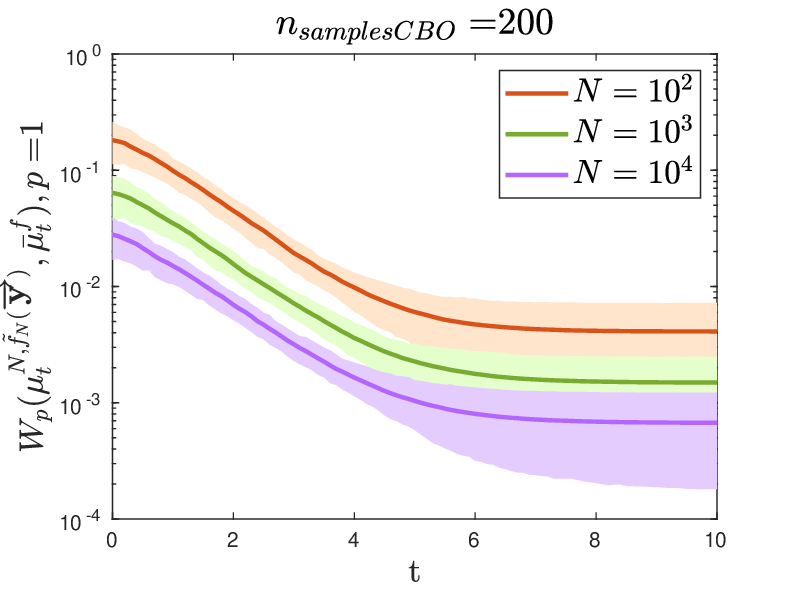}} 
	\subfigure[]{\includegraphics[scale = 0.49]{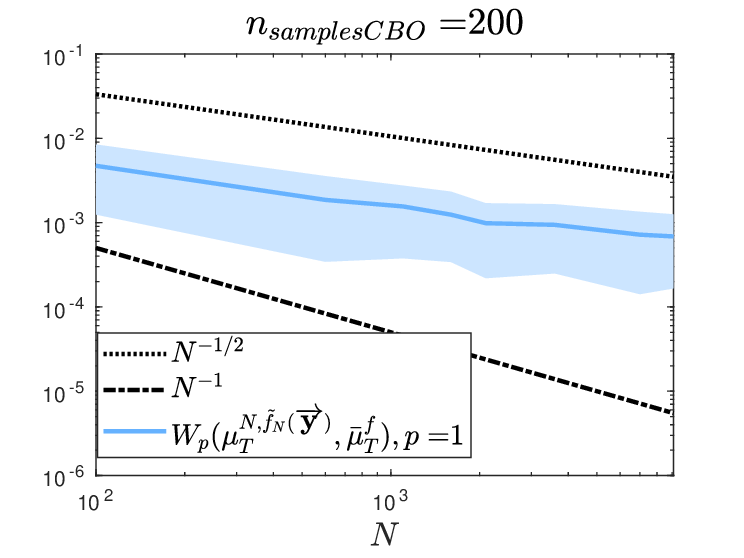}} 
	\caption{In the column on the left, evolution of the average error \eqref{def:newrate_quad} for $p=1$ (denoted by a continuous line) as a function of time for $N=10^2,10^3,10^4$. In the column on the right, evolution of the same error as a function of $N$ (varying from $10^2$ to $10^4$ with step $500$) and for $t=T$. The shaded area indicates the $[0.15,0.85]$ quantile interval related to the $n_{samplesCBO}$ runs of the CBO algorithm.
		Plot (a) has a log-scale on $y$-axis, plot (b) a log scale on both axis.
		The objectives considered are \eqref{def: ackleystocmodR1F} and \eqref{def: ackleystocmodR1f} with \eqref{param: ackleystocmodR1 Y}, parameters are set to \eqref{param5: test1,2}.}
	\label{fig5:ackleystocmodR1_quadcbo_rate2M_true_p1p2}
\end{figure}

We choose $F$ \eqref{def: ackleystocmodR1F} and $f$ \eqref{def: ackleystocmodR1f} and implement a $1$-dimensional isotropic CBO algorithm with parameters \eqref{param:tests1,2_complete} of Subsection \ref{subsecnum: rate SAA}.
We present our results in Figure \ref{fig5:ackleystocmodR1_quadcbo_rate2M_true_p1p2} for $p=1$.
\noindent Error \eqref{def:newrate_quad} yields a rate of $O(N^{-1/2})$ (see plot (b)). Denoting by $\bar{\mu}^{\tilde{f}_{N_{ref}}(\by)}_t$ the empirical measure associated to the sample $\{X^i_t(\by)\}^{i=1,\ldots,N_{ref}}$, solution to \eqref{eq3: complete cbo for ftildeN}, we may explain the rate by means of the following triangular inequality
\begin{equation}
	\label{eq:an_trueerror_quad}
	\underbrace{W_p \left(\mu^{N,\tilde{f}_N(\by)}_t, \bar{\mu}^f_t \right)}_{\tn{Joint-limit error \eqref{def:newrate_quad}}} \le 
	\underbrace{W_p \left(\mu^{N,\tilde{f}_N(\by)}_t, \bar{\mu}^{\tilde{f}_{N_{ref}}(\by)}_t \right)}_{\tn{Mean-field approximation error}} + 
	\underbrace{W_p \left(\bar{\mu}^{\tilde{f}_{N_{ref}}(\by)}_t, \bar{\mu}^f_t \right)}_{\substack{\text{Error independent of $N$,}\\ \text{dependent on $\tilde{f}_{N_{ref}}(\by) \approx f$}}}.
\end{equation}
The mean-field approximation error (first addend of the right-hand side of \eqref{eq:an_trueerror_quad}) has a rate $O(N^{-1/2})$ \cite{gerber2023mean}, while the second depends on the accuracy of the quadrature formula used, but not on $N$: all in all, the right-hand side of \eqref{eq:an_trueerror_quad} goes as $O(N^{-1/2})$, which is what we observe numerically \footnote{This remark justifies the choice of the composite midpoint rule rather than high order quadrature formulas.}.

\subsection{Test 3: Convergence rates for quadrature and CBO for $k \ge 2$}
\label{subsecnum: higher k}

In this section we investigate the impact of the dimension of the random space $k$ on error \eqref{def:newrate_quad}. 
It is well-known that a composite midpoint rule with $q^k=n$ nodes has an integration error of order $O(q^{-2/k}) = O(n^{-2/k^2})$ in the Euclidean distance provided that the integrand is $C^2$ in the $\mathbf{y}$-component, see e.g. \cite{hamming2012numerical}.
Then, in light of estimate \eqref{eq:an_trueerror_quad}, we expect the dependence of $k$ to influence only $W_p(\bar{\mu}^{\tilde{f}_{N_{ref}}(\by)}_t, \bar{\mu}^f_t)$: in particular, for $k>1$, we predict for \eqref{def:newrate_quad} to observe again a rate of $O(N^{-1/2})$, but with a bigger constant at the decay rate due to the lower accuracy of the quadrature formula for higher $k$.

We choose a class of objective functions where we can easily increase the dimension $k$ and which are inspired from the setting of linear-least squares problems \cite{hamming2012numerical}. We summarize our selection in Table \ref{tab:F_test3} for 
\begin{equation}
	\label{param:Flls}
	\rev{Y_1,Y_2,Y_3} \sim U \; \tn{independent, and $U \sim \mathcal{U}([0,2])$}.
\end{equation}

\begin{table}[h!]
	\centering
	\begin{tabular}{|c|c|c|}\hline
		$k$ & $F$ & $f$ \\ \hline
		$1$ & $F(x,Y_1) = (Y_1x)^2$ & $\EE[U^2] x^2$\\ \hline
		$2$ & $F(x,(Y_1,Y_2)^T) = (Y_1x -Y_2)^2$ & $\EE[U^2]x^2 -2(\EE[U])^2x+\EE[U^2]$ \\ \hline
		$3$ & 
		$ 
		F(x,(Y_1,Y_2,Y_3)^T) =  Y_1x^2 + Y_2x + Y_3
		$ & $\EE[U]x^2+\EE[U]x+\EE[U]$ \\ \hline
	\end{tabular}
	\caption{Objective functions for the case $3$. Please note that those do not fulfill the theoretical assumptions of Section \ref{sec:SAA}.}
	\label{tab:F_test3}
\end{table}

\begin{figure}[h!tbp]
	\centering
	\includegraphics[scale = 0.5]{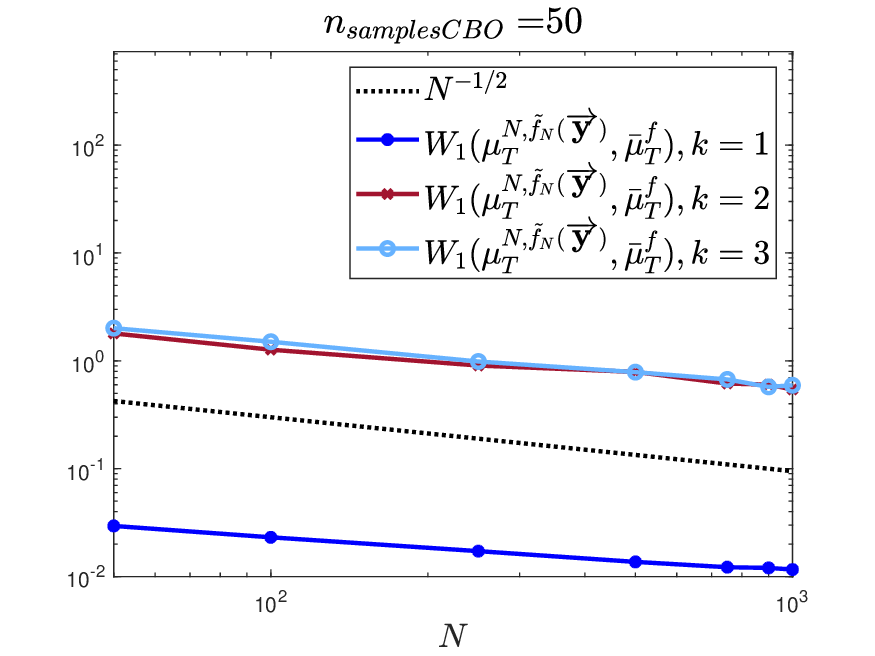}
	\caption{Evolution of average error \eqref{def:newrate_quad} as a function of $N$, for $p=1$, $t = T$.  The objectives considered are summarized in Table \ref{tab:F_test3} and the parameters are set as in \eqref{param:Flls} and \eqref{param5: test3}. 
		The plot has log scales on both axis.}
	\label{fig5:llsk123_quadcbo_comp_true_p1_d1bis}
\end{figure}

\noindent For this test, we use the isotropic CBO specified by the update rule \eqref{eq5: discretized cbo for M,N} with 
\begin{equation}
	\label{param5: test3}
	T = 7, \Delta t = 1; \quad \lambda = 1, \alpha = 40, \sigma = 0.5; \quad d=1, \mu_0 = \mathcal{U}([-3,3]^d).
\end{equation}
and we plot our results in Figure \ref{fig5:llsk123_quadcbo_comp_true_p1_d1bis} for $N_{ref}=10^3$, $N$ varying from $50$ to $1000$ and $p=1$. 
In the scenarios $k=1,2,3$ the same rate $O(N^{-1/2})$ is observed, with a constant in front of the decay rate that gets bigger as $k$ increases \footnote{The large magnitude of the error in Figure \ref{fig5:llsk123_quadcbo_comp_true_p1_d1bis} is due to the low number of iterations used, see \eqref{param5: test3}. This choice aims at speeding up the computational time as the evaluation of $x^{\alpha,\tilde{f}_N(\by)}_t$ requires $O(N^2) = O(N \cdot Q^k)$ computations for each iteration.}. 
This result reflects the conclusions drawn from estimate \eqref{eq:an_trueerror_quad} mentioned at the beginning of the section.

\subsection{Test 4: Success rates for a stochastic utility problem}
\label{subsecnum:test4_succrate}

In this section, we consider the stochastic utility problem proposed in \cite{nemirovski2009robust}.
We fix the dimension $k$ of the random space to be equal to the one of the search space $d$ and
\begin{equation}
	\label{def:Fcompl}
	F(x,(Y_1,\ldots,Y_d)^T) = \phi \left( \sum_{l=1}^d \left( \frac{l}{d} + Y_l \right) x_l \right) = \phi \left( x^T (a+\mathbf{Y})\right).
\end{equation}
$Y_1,\ldots,Y_d$ are independent normally distributed real-valued random variables with zero mean and unit variance, and $a \in \Rd$ is such that $a_l := \frac{l}{d}$, for any $l=1,\ldots,d$.
$\phi$ is a piecewise linear function given by
\begin{equation*}
	\phi(t) = \max \{ v_1+s_1t,\ldots,v_b+s_b t\}, \quad \tn{for $t \in \R$}
\end{equation*}
with $b \in \mathbb{N}$ and $v_1,\ldots,v_b,s_1,\ldots,s_b$ real-valued constants. In our experiments, we fix $b=4$ and
\begin{equation*}
	\mathbf{s} = (s_1,\ldots,s_4) = (-2,-1,1/2,1)^T; \quad \mathbf{v} = (v_1,\ldots,v_4) = (0,2,0,-1)^T,
\end{equation*}
so to have $b-1=3$ breakpoints $\mathbf{z} = (z_1,z_2,z_3) = (-2,4/3,2)^T$.
The advantage of the choice of $F$ \eqref{def:Fcompl} is that we can compute the value of the true objective $f(x)$ at a given candidate solution $x$ analytically and approximate the true minimizer $\xs$ of the corresponding optimization problem by bisection \cite{nemirovski2009robust}. 
In the following, we consider $k=d=1,2,3$ with the corresponding optimal points and values computed numerically, as stated above.
\begin{table}[h!]
	\centering
	\begin{tabular}{|c|c|c|}\hline
		$k=d$ & $\xs$ & $f(\xs)$ \\ \hline 
		$1$ & $0.82058 $ & $ 1.3927$\\ \hline
		$2$ & $ (0.35536,0.71572)^T$ & $1.3407$ \\ \hline
		$3$ & 
		$(0.20578,0.40601,0.61735)^T $ & $1.2895$ \\ \hline
	\end{tabular}
	\caption{$\xs$ and $f(\xs)$ for $k=d=1,2,3$ and $F$ given by \eqref{def:Fcompl}.}
	\label{tab:xsfxs_test4}
\end{table}

As mentioned in Subsection \ref{subsecnum: higher k}, the rate of any deterministic quadrature formula depends on the regularity of the integrand function and on the dimension of the random space $k$; in contrast, the exponent $1/2$ of the SAA error is independent of $k$ (equal to $O(M^{-1/2})$: robust, but in general slower) and on the smoothness of $F$ in the $\mathbf{y}$-component, see \cite{shapiro2003monte,shapiro2021lectures,caflisch1998monte}.
In this section, we compare the performances of algorithms \eqref{eq2: complete cbo for fhatM} (SAA and CBO) and \eqref{eq3: complete cbo for ftildeN} (quadrature and CBO) and aim at numerically testing the just-mentioned theoretical result.
To carry out the comparison, we evaluate the success rate at the final time $T$ for algorithms \eqref{eq2: complete cbo for fhatM} (SAA and CBO) and \eqref{eq3: complete cbo for ftildeN} (quadrature and CBO). 
More precisely, in agreement with \cite{pinnau2017consensus,carrillo2021consensus} and given a threshold $thr >0$, we consider a run successful for algorithm \eqref{eq2: complete cbo for fhatM} if the candidate minimizer of the SAA and CBO, namely $\EE_{\PP}[x^{\alpha,\hat{f}_M(\bY(\cdot))}_T]$, is contained in the open $||\cdot||_{\infty}$-ball with radius $thr$ around the true minimizer $\xs$ (hereafter denoted by $B_{thr}^{\infty}(\xs)$): as in Subsection \ref{subsecnum: rate SAA}, we added the expectation $\EE_{\PP}$ to accommodate for the uncertainty of $\bY(\cdot)$ and numerically simulate it by averaging over $n_{samplesY}=100$ realization. Then, we consider a run successful for algorithm \eqref{eq3: complete cbo for ftildeN} if $x^{\alpha,\tilde{f}_N(\by)}_T \in B_{thr}^{\infty}(\xs)$. We compute the success rates with $n_{samplesCBO} = 100$ realizations of the CBO algorithms.\\

We use an anisotropic \footnote{The choice of anisotropic random exploration process for problems with a high dimensional search space has been shown to be more competitive than the isotropic one, thanks to the independence of the parameter constraints of the dimensionality $d$ \cite{carrillo2021consensus,fornasier2022anisotropic}.} CBO algorithm with parameters \eqref{param5: test1,2} and thresholds $thr_1 = 0.50, thr_2 = 0.25, thr_3 = 0.10$. We require $N=M=Q^k$: the equal choice of number of agents, samples and nodes allows for a comparison between the SAA and quadrature approaches. 
We summarize our results in Tables \ref{tab:test4_saacbo} and \ref{tab:test4_quadcbo} for $N \in [100,500,1000]$ \footnote{We numerically assessed that $N=50$ was the minimum value to get significant success rates.}. The cells of the table have been colored in white for a success rate of $100\%$, green for the range $90-99\%$, yellow for the range $50-89\%$, and red for any percentage below $49\%$.

\begin{table}[h!]
	\centering
		\begin{tabular}{|c|c|c|c|c|c|c|c|c|c|}\hline
			& \multicolumn{3}{c|}{$k=d=1$} & \multicolumn{3}{c|}{$k=d=2$} &  \multicolumn{3}{c|}{$k=d=3$}  \\ \hline 
			$N$ & $thr_1$ & $thr_2$ & $thr_3$ & $thr_1$ & $thr_2$ & $thr_3$ & $thr_1$ & $thr_2$ & $thr_3$ \\ \hline
			$100$ & $100\%$ & $100\%$ & $100\%$ & $100\%$& $100\%$ &  $100\%$& $100\%$ & $100\%$ & \cellcolor[HTML]{B8FF33} $99\%$ \\ \hline
			$500$ & $100\%$ & $100\%$ & $100\%$ & $100\%$ & $100\%$& $100\%$ & $100\%$ & $100\%$ & $100\%$ \\ \hline
			$1000$ & $100\%$ & $100\%$ & $100\%$ & $100\%$ & $100\%$& $100\%$ & $100\%$ & $100\%$ & $100\%$ \\ \hline
		\end{tabular}
	\caption{Success rate at final time $T$ for algorithm \eqref{eq2: complete cbo for fhatM} applied to function $F$ \eqref{def:Fcompl} and parameters \eqref{param5: test1,2}.}
	\label{tab:test4_saacbo}
\end{table}
\begin{table}[h!]
	\centering
	\begin{tabular}{|c|c|c|c|c|c|c|c|c|c|}\hline
		& \multicolumn{3}{c|}{$k=d=1$} & \multicolumn{3}{c|}{$k=d=2$} &  \multicolumn{3}{c|}{$k=d=3$}  \\ \hline 
		$N$ & $thr_1$ & $thr_2$ & $thr_3$ & $thr_1$ & $thr_2$ & $thr_3$ & $thr_1$ & $thr_2$ & $thr_3$ \\ \hline
		$100$ & $100\%$ & $100\%$ & $100\%$ & \cellcolor[HTML]{B8FF33} $99\%$ & \cellcolor[HTML]{FF7D33} $0\%$ & \cellcolor[HTML]{FF7D33} $0\%$ & \cellcolor[HTML]{FF7D33}$23\%$ & \cellcolor[HTML]{FF7D33}$2\%$ & \cellcolor[HTML]{FF7D33}$0\%$ \\ \hline
		$500$ & $100\%$ & $100\%$ & $100\%$ & $100\%$ & $100\%$ & \cellcolor[HTML]{FF7D33}$0\%$ & \cellcolor[HTML]{FFEC33}$64\%$ & \cellcolor[HTML]{FF7D33}$0\%$ & \cellcolor[HTML]{FF7D33}$0\%$\\ \hline
		$1000$ & $100\%$ & $100\%$ & $100\%$ & $100\%$ & $100\%$ & $100\%$ & \cellcolor[HTML]{FFEC33}$75\%$ & \cellcolor[HTML]{FF7D33}$0\%$ & \cellcolor[HTML]{FF7D33}$0\%$ \\ \hline
	\end{tabular}
	\caption{Success rate at final time $T$ for algorithm \eqref{eq3: complete cbo for ftildeN} applied to the function $F$ \eqref{def:Fcompl} and parameters \eqref{param5: test1,2}.}
	\label{tab:test4_quadcbo}
\end{table}

\noindent In both tables, we assess an increase of the success rate for increasing threshold values and $N=M=Q^k$.
Having a closer look at rows $3$ to $5$ of Table \ref{tab:test4_saacbo}, we observe a success rate of $100\%$ for all $k=d$: this is consistent with the theoretical result mentioned in the previous page, that stated that the SAA rate is $O(M^{-1/2})$, independent of $k$.
On the other hand, if we shift our attention to rows $3$ to $5$ of Table \ref{tab:test4_quadcbo}, we see that a higher $N$ is required to obtain a $100\%$ success rate and that such rate is highly influenced by the dimension of the random space $k$: again, this is consistent with the result that the accuracy of the composite midpoint quadrature formula (that we are using in the simulations) depends on $k$ and on the regularity of $F$. 

\section{Summary and conclusion}
\label{sec: saa+cbo vs quad+cbo}

In this work we propose two approaches for solving an optimization problem where the cost function is given in the form of an expectation by means of a stochastic particle dynamics based on consensus. 
The two methods share the property of replacing the objective function with a suitable approximation, which we chose to be a Monte Carlo type estimator based on a fixed sample and a quadrature formula, respectively. 
In light of the numerical results of Section \ref{sec:numerics} and of \cite{albi2021trails, shapiro2003monte, shapiro2021lectures,pareschi2013interacting}, in which comparisons between Monte Carlo (MC) type and deterministic approaches are carried out, we now proceed to summarize the main advantages and disadvantages of the two approaches presented in the manuscript.

The MC method may be applied even if the law $\nuy$ of the random vector $\mathbf{Y}$ is unknown or, when $\nuy$ is absolutely continuous with respect to the Lebesgue measure and it admits the density $\theta_\mathbf{Y}$, if $\theta_\mathbf{Y}$ is non-smooth; on the other hand, the quadrature approach requires the existence of $\theta_\mathbf{Y}$ and its smoothness. 

As regards the computational complexity, we start by observing that in the stochastic approach we have a loop on $n_{samplesY}$ to accommodate for the stochasticity of the random variable $\mathbf{Y}$, which disappears in presence of the deterministic approach. The cost of each iteration of the CBO algorithm is reported in Table \ref{tab6: comp complexity}.

\begin{table}[h!]
	\centering
	\begin{tabular}{|c|c|c|}\hline
		& \makecell{SAA approach,\\ algorithm \eqref{eq2: complete cbo for fhatM}} & \makecell{Quadrature approach,\\ algorithm \eqref{eq3: complete cbo for ftildeN}} \\ \hline
		\makecell{ Number of evaluations of $F$\\ to compute $\hat{f}_M,\tilde{f}_N$}  & \makecell{$O(N \cdot M)$\\ to compute $\hat{f}_M$} & \makecell{$O(N \cdot Q^k)$ \tablefootnote{We observe that we may reduce such complexity to $O(N \cdot Q(\log Q)^{k-1})$ by using sparse grids \cite{bungartz2004sparse}.}\\ to compute $\tilde{f}_N$} \\ \hline
		Number of evaluations of $\theta_\mathbf{Y}$ & $-$ & $O(Q^k)$ \\ \hline
		Computation of $x^{\alpha,\hat{f}_M}_{t_h}, x^{\alpha,\tilde{f}_N}_{t_h}$ & \multicolumn{2}{c|}{$O(N)$}  \\ \hline
	\end{tabular}
	\caption{Analysis of the computational complexity for each iteration $h$ of the CBO algorithm for the SAA and quadrature approaches with update rule given by \eqref{eq5: discretized cbo for M,N}.}
	\label{tab6: comp complexity}
\end{table}

\noindent As mentioned in Subsections \ref{subsecnum: higher k} and \ref{subsecnum:test4_succrate} and as numerically checked in Tables \ref{tab:test4_saacbo} and \ref{tab:test4_quadcbo} of Subsection \ref{subsecnum:test4_succrate}, the regularity of the integrand function and the dimension of the random space $k$ have an impact on the rate of quadrature formulas, while the exponent $1/2$ of the MC error is independent of them. Then, merging this information with the analysis carried out in Table \ref{tab6: comp complexity}, we conclude that the quadrature approach is preferred in the scenario of low dimensional random space ($k=1$), as, for a comparable success rate, there isn't the need to have a loop on $n_{samplesY}$ and the number of evaluations of $F$ is comparable with the ones of the SAA approach. However, as $k$ increases, one may want to opt for the SAA approach where the number of evaluations of $F$ is independent of $k$.

In this work we also conduct a theoretical analysis at the mean-field level and prove that the long time behavior of the SAA approach coincides with the one of the true problem. We assess that the identical choice of number of particles and nodes in the quadrature method leads to a mean-field equation in the extended phase space.
We define two metrics for the former approach and one for the latter and investigate the rates of convergence of the methods: we conclude that they are consistent with the results of the classical theory of SAA and quadrature formulas.\\

We conclude the section by mentioning variable-sample strategies \cite{homem2003variable}, which could be a third approach to solve stochastic programming problem \eqref{eqi: min problem main}.
In contrast to the SAA approach, where a sample $\by$ is fixed at the beginning and then a minimization problem for $\hat{f}_M(\by)$ is solved, in variable-sample strategies a sample is drawn at each iteration of the algorithm used to solve the minimization problem, so that $\hat{f}_M$ depends on the iterate $h$.
It is still open how to include this strategy in the previous work.


\section*{Acknowledgments}
The authors thank the Deutsche Forschungsgemeinschaft (DFG, German Research Foundation) for the financial support through $442047500/$ SFB1481 within the projects B04 (Sparsity fördernde Muster in kinetischen Hierarchien), B05 (Sparsifizierung zeitabhängiger Netzwerkflußprobleme mittels diskreter Optimierung) and B06 (Kinetische Theorie trifft algebraische Systemtheorie), for the financial support under Germany’s Excellence Strategy EXC-2023 Internet of Production 390621612 and under the Excellence Strategy of the Federal Government and the Länder and for support received funding from the European Union’s Horizon Europe research and innovation programme under the Marie Sklodowska-Curie Doctoral Network Datahyking (Grant No. 101072546).
The work of SB is funded by the Deutsche Forschungsgemeinschaft (DFG, German Research Foundation) – 320021702/GRK2326 – Energy, Entropy, and Dissipative Dynamics (EDDy).

\bibliographystyle{plain}
\bibliography{references_arXiv_R3}

\end{document}